\documentclass[11pt]{article}
\usepackage[title]{appendix}

\usepackage{authblk}
\usepackage[utf8]{inputenc}
\usepackage{amsmath,amssymb,amsthm,amsfonts,tikz,gensymb, phonetic, float, caption, subcaption, comment, mathtools,pgfplots}
\usepackage{xurl}
\usepackage{pifont}
\usepackage{exscale}
\usepackage{graphicx}
\usepackage[font=footnotesize]{caption}
\usepackage{subcaption}
\usepackage[margin=1in]{geometry}
\usepackage{hyperref}
\usepackage{enumitem}
\usepackage{siunitx}
\usepackage{parskip}
\usepackage{listings}
\usepackage{color,pxfonts,fix-cm}
\usepackage{latexsym}
\usepackage[mathletters]{ucs}
\DeclareUnicodeCharacter{32}{$\ $}
\DeclareUnicodeCharacter{46}{\textperiodcentered}

\usepackage{textcomp, gensymb}

\usepackage[compact]{titlesec}        
\titlespacing{\section}{0pt}{12pt}{0pt}
\AtBeginDocument{
  \setlength\abovedisplayskip{0pt}
  \setlength\belowdisplayskip{0pt}}
\providecommand{\keywords}[1]{\textbf{Keywords:} #1}

\newcommand{\T}{\mathsf{T}}

\newcommand{\dx}{\Delta x}
\newcommand{\dy}{\Delta y}
\newcommand{\dz}{\Delta z}

\newcommand{\eq}[1]{Equation~#1}
\newcommand{\tabRef}[1]{Table~#1}
\newcommand{\secRef}[1]{Section~#1}
\newcommand{\figRef}[1]{Figure~#1}
\newcommand{\ppt}{\mathrm{ppt}}
\renewcommand{\Vec}[1]{\mathbf{#1}}

\title{Investigating Ocean Circulation Dynamics Through Data Assimilation: A Mathematical Study Using the Stommel Box Model with Rapid Oscillatory Forcings}

\author[1]{Nathaniel Smith}
\author[2]{Anvaya Shiney-Ajay}
\author[3]{Emmanuel Fleurantin}
\author[4]{Ivo Pasmans}
\date{}

\affil[1]{\small{Miami University, Oxford, OH 45056}}
\affil[2]{\small{University of Minnesota, Minneapolis, MN 55455}}
\affil[3]{\small{George Mason University, Fairfax, VA 22030}}
\affil[4]{\small{University of Reading, Reading RG6 6AH, United Kingdom}}
\affil[*]{\footnotesize{email: smithni@miamioh.edu, ajayx006@umn.edu, efleuran@unc.edu, i.c.pasmans@reading.ac.uk}\vspace{-0.5em}}

\pgfplotsset{compat=1.18}

\begin{document}

\maketitle

\begin{abstract}
We investigate ocean circulation changes through the lens of data assimilation using a reduced-order model. Our primary interest lies in the Stommel box model which reveals itself to be one of the most practicable models that has the ability of reproducing the meridional overturning circulation. The Stommel box model has at most two regimes: TH (temperature driven circulation with sinking near the north pole) and SA (salinity driven with sinking near the equator). Currently, the meridional overturning is in the TH regime. Using box-averaged Met Office EN4 ocean temperature and salinity data, our goal is to provide a probability that a future regime change occurs and establish how this probability depends on the uncertainties in initial conditions, parameters and forcings. We will achieve this using data assimilation tools and DAPPER within the Stommel box model with fast oscillatory regimes. 
\end{abstract}

\begin{flushleft}
\\

\vspace{1em}

\end{flushleft}

\section{Introduction}\label{sec:intro}
Thermohaline circulation is the global current driven by density variations in the oceans. Its branch in the Atlantic ocean is known as the Atlantic Meridional Overturning Circulation (AMOC). In the northern part of the Atlantic it transports warm, salty water from the subtropics to the Arctic. Here it comes into contact with colder (sub)polar waters and sink forming the North-Atlantic Deep Water and flows southward in the deep ocean \cite{talley_descriptive_2011}. It plays a crucial role in the climate system, especially for the climate in northwestern Europe. Based on a simple two-box model, \cite{stommel} suggested that AMOC can be described as a bistable system characterized by two distinct equilibrium states: TH (thermally-driven), representing the current circulation regime with sinking near the pole, and SA (salinity-driven), signifying a reversed circulation with equatorial sinking. The existence of multiple equilibria has since been confirmed by both reversal visible in the paleo-record as well as more sophisticated climate models, see e.g. \cite{bryan_high-latitude_1986,marotzke_instability_1988,maier-reimer_experiments_1989}. 

In the region where sinking occurs, temperature and salinity exert opposing effects on density. On one hand, heat exchange with the cold (sub)polar environment cools down the subtropical water that are being advected northward, increasing its density. On the other, contact with melt water freshens it, decreasing its density. Currently, the former effect dominates. However, it is feared that climate change driven increases in temperature, ice melt and precipitation can tip the balance to the latter, thus reversing the circulation. The transition between the circulation regimes can occur through various mechanisms \cite{doi:10.1073/pnas.0802430105,https://doi.org/10.1029/2010GL044486,scheffer_early-warning_2009}. Beyond a critical point, parameter values may lead to bifurcation-induced tipping (B-tipping). Alternatively, the rate of parameter changes might trigger rate-induced tipping (R-tipping). Furthermore, system variability due to inherent noise could result in noise-induced tipping (N-tipping). Previous research, as documented in \cite{weijer_stability_2019, RegimeTrans, Otto2354, cimatoribus_meridional_2014, Cast, MR2276274, MR4165514, MR2664454}, has established thresholds for tipping occurrence. In \cite{esd-12-819-2021}, the Stommel box model is coupled to a sea-ice component in the polar box, and displays rate-induced tipping. 

It is still an open question how close to such a tipping point the world is \cite{weijer_stability_2019}. In 2023 alone, two studies on the topic \cite{ditlevsen_warning_2023,van_westen_physics-based_2024} have received considerable attention in the press, but their results have been criticised \cite{ben-yami_uncertainties_2023}. The aim of this study is twofold. First, to explore whether it is possible to use observations of ocean temperature and salinity in combination with a reduced-order model in a process called data assimilation to estimate model parameters. Second, determine how the addition of forcing, both seasonal and climatological, to the Stommel model affects this probability of an AMOC reversal. 

It is important to note that while one might assume that a two-box model would have limitations in fully describing tipping behavior, our approach, incorporating high-frequency oscillations, parameter variability, and data assimilation, brings us close to realistic scenarios. This methodology also sheds light on the feasibility of assimilating data to further understand the qualitative dynamics of our model. It offers a window of validity for our predictions, considering factors like the strength of the forcing.

{\subsection{Leveraging Conceptual Models for Parameter Estimation and Predictive Analysis} The use of conceptual models, such as the Stommel box model, for parameter estimation and prediction of complex systems like the AMOC is motivated by several important considerations. Firstly, reduced-complexity models often provide clearer insights into the fundamental mechanisms driving the system's behavior, which can be obscured in more complex models \cite{sherman2023data,smith2013data,sevellec2015unstable}. This interpretability is vital for understanding the key processes governing AMOC dynamics. 

In light of these advantages, our objective is to model the AMOC and examine its behavior under various conditions. We aim to assess the feasibility of estimating temperature diffusivity, salinity diffusivity, and advection parameters by assimilating ocean temperature and salinity data. Additionally, we seek to forecast the probability of circulation regime changes. The simplicity of the chosen conceptual model makes this analysis possible as we are able to initialize large ensembles at low computational cost. It also serves as a baseline in complexity with more detailed versions of the model available for further studies on higher order dynamics \cite{hawkins2011bistability, jackson2023understanding}.

Lastly, these models serve as useful benchmarks against which more complex models can be evaluated, helping to assess whether increased complexity leads to improved predictive skill \cite{variabilityDG,chapman_tipping_2024,cimatoribus_meridional_2014,boers2021observation,dekker2018cascading}. While we acknowledge the limitations of conceptual models, their use in conjunction with data assimilation techniques offers a powerful approach to understanding AMOC dynamics and predicting potential tipping points.}


\subsection{Outline of the Paper}

The outline of the paper is as follows. First, we introduce the autonomous and non-autonomous Stommel model in Section \ref{mat_stom} and study the qualitative behavior of the solutions as time goes to infinity. In Section \ref{sec:dap}, we explain our data assimilation method and the assimilated observations. In Section \ref{sec:twin}, we conduct a twin experiment to investigate the potential of data assimilation to estimate model parameters. In Section \ref{sec:obs}, observations based on the Met Office Hadley Centre EN4 data set are assimilated in a non-autonomous Stommel model and the impact of the different forcings on the tipping probability is investigated. In Section \ref{sec:future}, we present and analyze the results of our data assimilation experiments, examining the model's behavior under various climate change scenarios and assessing the likelihood of AMOC regime shifts. Finally in Section \ref{sec:concl}, we discuss our findings. 

\section{Mathematical Analysis of the Stommel box Model: Autonomous and Non-Autonomous Perspectives}\label{mat_stom}

\subsection{The Stommel Box Model}
\label{sec:stommel_model}

The Stommel model comprises two fully mixed water containers. One represents the region of the North-Atlantic just north of the equator with temperature and salinity denoted as $(T_{e*}, S_{e*})$, while the other symbolizes the polar region with temperature and salinity as $(T_{p*}, S_{p*})$. (The asterisk subscripts here denote variables with a unit.) These containers are interconnected by both a capillary tube and an overflow area, facilitating water circulation while preserving a constant volume in each container. The relative salinities and temperatures of each box determine the flow of water, and when the flow of water reverses, we consider that to be a regime change.

\begin{figure}[h!]
\centering
\includegraphics[scale=0.25]
{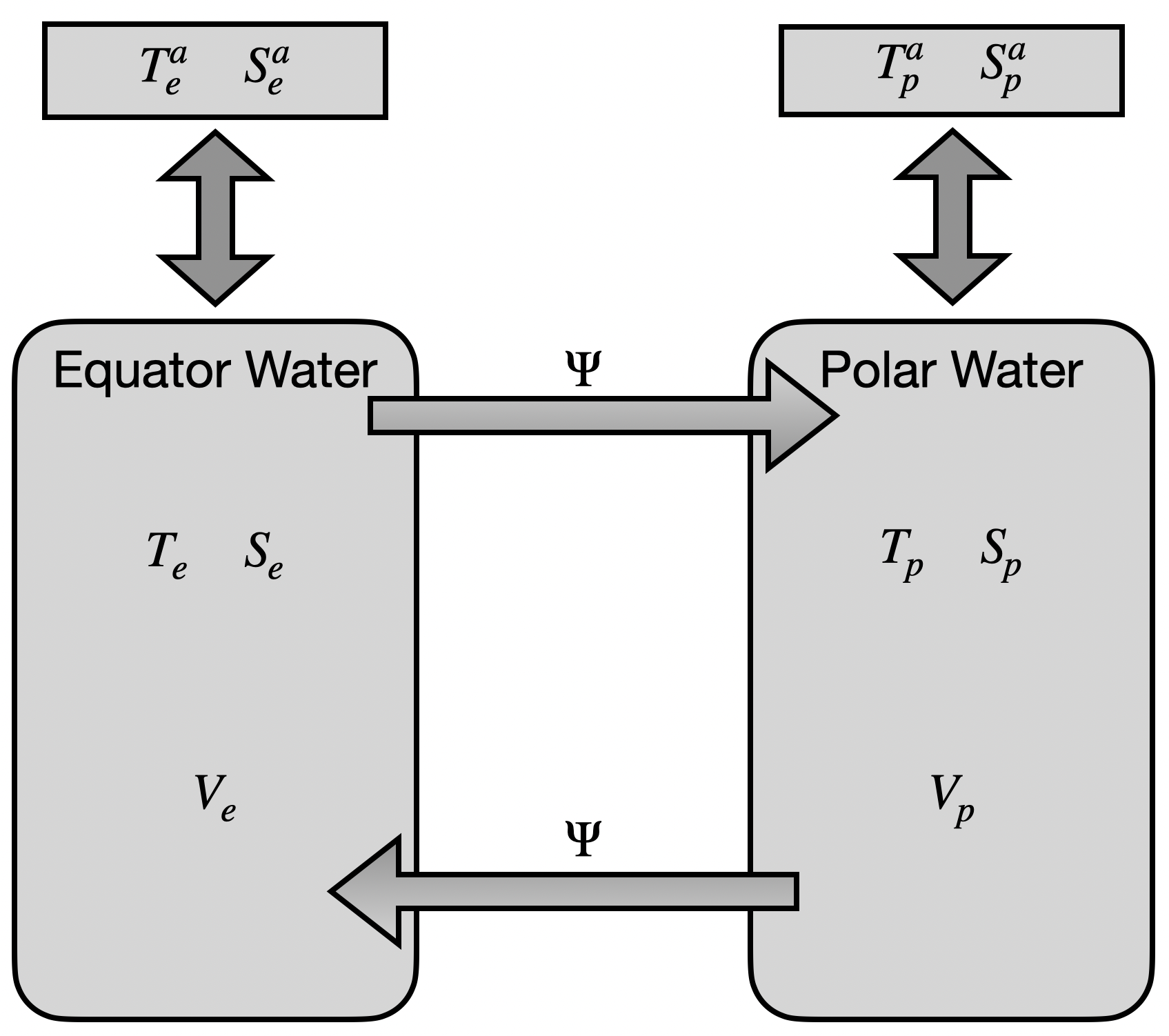}
\caption{\raggedright A conceptual depiction of how Henry Stommel modeled thermohaline circulation. The circulation between the two connected bodies of well-mixed water in the Stommel model is driven by density gradients between the two boxes. }
\end{figure}

Initially, then, this scenario is described with four equations which are listed below: 

\begin{equation}
\begin{split}
    \dx \dz \dy_{e} \frac{dT_{e*}}{dt_*} = & k_{T} \dx \dy_{e} (T_e^a - T_{e*}) + |\Psi_*|(T_{p*} - T_{e*}),\\ 
    \dx \dz \dy_{e} \frac{dS_{e*}}{dt_*} = & k_{S} \dx \dy_{e} (S_e^a - S_{e*}) + |\Psi_*|(S_{p*} - S_{e*}), \\ 
    \dx \dz \dy_{p} \frac{dT_{p*}}{dt_*} = & k_{T} \dx \dy_{p} (T_p^a - T_{p*}) + |\Psi_*|(T_{e*} - T_{p*}),\\
    \dx \dz \dy_{p} \frac{dS_{p*}}{dt_*} = & k_{S} \dx \dy_{p} (S_p^a - S_{p*}) + |\Psi_*|(S_{e*} - S_{p*}) - m^{*} S_{p*},
\end{split}
\label{eq:main}
\end{equation}

with $\dz$ and $\dx$ being the depth and the zonal width of the boxes, $\dy_{p}$ and $\dy_{e}$ the meridional width of the  polar and equatorial box respectively, $k_{T}$ the temperature diffusion coefficient, and $k_{S}$ the salinity diffusion coefficient. The model incorporates the exchange of heat and salt within each box as a response to surface forcing, employing relaxation toward specified ocean surface temperature and salinity values $(T^a, S^a)$. Next to this, the model accomodates for the addition of fresh water supply to the polar box at a discharge rate of $m^{*}$. {In particular, this quantity represents the melting of icebergs due to climate change and the release of additional fresh water that comes from that.}

The flow rate is given by the following equation:
\begin{equation} 
    \Psi_* = \gamma \dx \dz (\rho_{p*} - \rho_{e*})/\rho_0,
    \label{eq:trans}
\end{equation}
where $\gamma$ is the advection coefficient and $\rho_0$ is a reference density. Since those are both constants, the flow rate is linearly related to the difference in water density between the two boxes. We note that if the flow is due to denser water in the polar box, as it is for the current TH regime, the flow rate is taken to be positive. 

The Stommel system also assumes the following linear equation of state 
\begin{equation} 
    \rho_* = \rho_0(1 - \alpha_T(T_* - T_0) + \alpha_S(S_* - S_0)),
    \label{eq:eos}
\end{equation}
where a ``0" subscript indicates reference value, $\alpha_T$ is the thermal expansion coefficient and $\alpha_S$ the haline coefficient. 

But by instead looking only at the temperature and salinity relative to one another, we can reduce this down to a two variable system by subtracting the equations for the pole box from the equatorial box in the case that $m^{*}=0$:
\begin{equation}\label{eq:dimensionless}
\begin{split}
\frac{dT}{dt} =& \eta_1 - T(1+|T-S|),\\
\frac{dS}{dt} =& \eta_2 - S(\eta_3+|T-S|),\\
\Psi =& T - S,
\end{split}
\end{equation}
where 
\begin{equation*}
\begin{aligned}
  &\begin{aligned}
    T &= \frac{\dz \gamma \alpha_{T}}{\overline{\dy} k_{T}} (T_{e*}-T_{p*}) \\[1ex]
    S &= \frac{\dz \gamma \alpha_{S}}{\overline{\dy} k_{T}}(S_{e*}-S_{p*}) \\[1ex]
    \Psi &= \frac{1}{\overline{\dy} \dx k_{T}} \Psi_{*}
  \end{aligned}
  \qquad
  &\begin{aligned}
    t &= \frac{k_{T}}{\Delta z} t_{* } \\[2.5ex]
    \overline{\dy} &:= \frac{\dy_{p} \dy_{e}}{\dy_{p}+\dy_{e}}
  \end{aligned}
\end{aligned}
\end{equation*}
.

Here, $\eta_1,\eta_2$, and $\eta_3$ are parameters independent of the ocean state. Their full values are as follow:
\begin{align}
    \eta_1 &= \frac{\dz \gamma \alpha_{T}}{\overline{\dy} k_{T}} (T^a_e-T^a_p) \label{eq:eta1}\\
    \eta_2 &= \frac{k_S}{k_T} \frac{\dz \gamma \alpha_{S}}{\overline{\dy} k_{T}}(S^a_e-S^a_p) \label{eq:eta2}\\
    \eta_3 &= \frac{k_S}{k_T}
    \label{eq:eta3} 
\end{align}

Setting $\frac{dT}{dt}=0$, $\frac{dS}{dt}=0$ in \eq{\eqref{eq:dimensionless}}, we get the following equilibrium solutions: 
\begin{align*}
    && T = \frac{\eta_1}{1 + |\Psi|} && \text{and} && S = \frac{\eta_2}{\eta_3 + |\Psi|} \,, &&
\end{align*}
In the TH regime, we have $\Psi=|\Psi|$. In this case, the effect of the temperature difference dominates the effect of the salinity difference and the system is thermally driven. We find that 
\begin{equation}\label{eq:Stom1}
 \eta_2 = -\Psi^2 - \eta_3\Psi + \eta_1\left(\frac{\eta_3 + \Psi}{1 + \Psi}\right) \,.
\end{equation} 
On the other hand, in the SA regime, we have $\Psi=-|\Psi|$. In this case, the transport due to the salinity difference dominates over the transport due to temperature difference. We find that 
\begin{equation}\label{eq:Stom2}
    \eta_2 = \Psi^2 - \eta_3\Psi + \eta_1\left(\frac{\eta_3 - \Psi}{1 - \Psi}\right)\,.
\end{equation}

The solutions from \eq{\eqref{eq:Stom1}} and \eq{\eqref{eq:Stom2}} for fixed $\eta_1$, $\eta_3$ can be shown together as a function of the dimensionless transport. From this, it is clear that for a range of $\eta_2$, three steady-state solutions for the transport are possible {, since the vertices of both parabolas described above lie in the area where $\psi>0$, meaning the TH regime provides two solutions for certain values of $\eta_2$ and the SA regime provides one.}

\subsection{The Non-Autonomous Stommel Model}

We extend the Stommel model to include seasonal variations, represented by periodic changes in surface temperatures and salinities. {We also know that any oscillatory behavior will not arise from the autonomous Stommel model (see Appendix A)} and thus will solely emerge due to the fact that $\eta_{1}$ and $\eta_{2}$ are no longer constant in time but oscillate, i.e. 
\begin{align}
\eta_1 =& \eta_{1,auto} + B\sin(\Omega t) \label{eq:eta1_B},\\
\eta_2 =& \eta_{2,auto} + \hat{B}\sin(\Omega t) ,\label{eq:eta2_Bhat}
\end{align}
with $\eta_{\cdot,auto}$ representing the values used in the autonomous system in section~\ref{eq:dimensionless}. 
Here, $B$ and $\hat{B}$ are distinct, but we define $A=B-\hat{B}$ for future use. Using that, our definition of $\Psi$, and the overall Stommel system, we can write:
\[\begin{split}
\frac{d}{dt}\Psi=&\frac{d}{dt}(T-S)=\frac{dT}{dt}-\frac{dS}{dt}\\
=&(\eta_{1,auto} + B\sin(\Omega t) - T(1+|T-S|))
-(\eta_{2,auto} + \hat{B}\sin(\Omega t) - S(\eta_3+|T-S|))\\
=&\eta_{1,auto} - \eta_{2,auto} + \eta_3 S - T - (T-S)|T-S| 
+ (B-\hat{B})\sin(\Omega t)\\
=&\eta_{1,auto} - \eta_{2,auto} + \eta_3 S - T - \Psi|\Psi| + A\sin(\Omega t),
\end{split}\]
then we substitute in $S=T-\Psi$, which ultimately gets us:
\begin{equation}\label{eq:nonaut}
    \begin{split}
        \frac{d}{dt}{\Psi}=&\eta_{1,auto} - \eta_{2,auto} + (T-\Psi)\eta_3 - T - \Psi|\Psi| + A\sin(\Omega t),\\
        \frac{d}{dt}{T}=&\eta_{1,auto} + B\sin(\Omega t) - T(1+|\Psi|).
    \end{split}
\end{equation}
for a system without climate change, i.e., a system in which $m^{*}=0$. 

\begin{figure}[h!]
\centering
\includegraphics[scale=0.5]{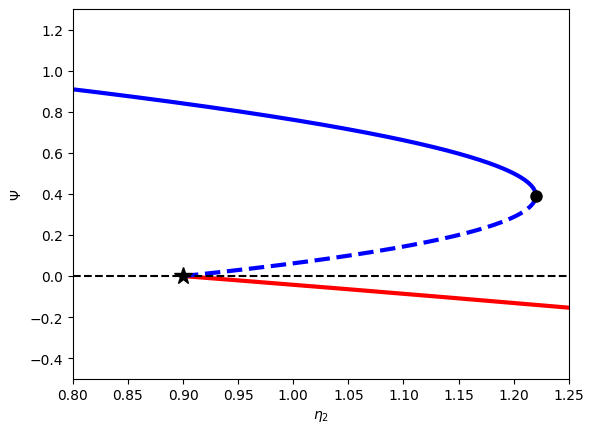}
\caption{\raggedright A $\Psi$-vs.-$\eta_{2,auto}$ bifurcation diagram for the system in equation \ref{eq:nonaut} with $\eta_{1,auto}=3$, $\eta_3=0.1$, and A=B=0. The solid red and blue lines are stable equilibria, while the dotted blue line represents an unstable equilibrium. The solid black circle indicates the system's saddle-node bifurcation, while the black star indicates the NSF bifurcation. The dotted black line emphasizes the discontinuity set $\Sigma_0$.
}
\label{fig:bif1}
\end{figure}

In order to get a full picture of the dynamics of \eq{\eqref{eq:nonaut}}, it is thus fitting to analyze its $\Psi$-vs.-$\eta_{2,auto}$ bifurcation diagram. It's worth noting that the discontinuity surface $\Sigma_0 =  \{(T,S) \in \mathbb{R}^2_+ | T=S\}$ give rise to a Non-Smooth Fold (NSF) bifurcation (marked by a black star in \figRef{\ref{fig:bif1}}), where the stable equilibrium branch terminates at $\eta_{2c}\equiv \eta_1\eta_3$. 

{Understanding NSF bifurcations can be helpful in predicting critical transitions or tipping points in complex systems, aiding in risk assessment and management. This paper examines the non-autonomous Stommel box model as a foundation for studying transitions in thermohaline circulation, while providing insight into the dynamical impact of NSF bifurcations. 

In the context of our model,} the NSF arises when the unstable equilibrium solution (dotted blue curve) and the focus (solid red curve) intersect with the discontinuity surface. This tipping point corresponds to the rapid transition from solutions near the temperature-dominated branch of nodes (solid blue curve) to the salinity-dominated branch of focus points{. }

For the non-autonomous Stommel model (\eq{\eqref{eq:nonaut}}), {under high frequency forcing, we observe amplification or attenuation of patterns in $\psi$ which depends linearly on the amplitude-to-frequency ratio $(A/\Omega)$.} 
{That is to say, solutions starting near the NSF bifurcation point may naturally shift away towards higher values of $\eta_2$, or be driven towards lower values under the effect of data assimilation.} The foundational work for this analysis can be found in a paper by Budd, Griffith, and Kuske \cite{budd_dynamic_2022}. Their study provides valuable insights into the expected behaviors of the parameters $\eta_1$ and $\eta_2$, as well as the equilibria, when incorporating seasonal variations. Building upon their findings, we extend the investigation by introducing data assimilation techniques to further enhance our understanding of these complex dynamics. By assimilating observational data into the model, we aim to explore how the combination of seasonal variations and data assimilation influences the system's behavior and equilibria.

\section{Inversion of the Stommel box model}\label{sec:dap}

The solution of the Stommel box model in \ref{eq:main} depends on the surface diffusion parameters $k_{T}$ and $k_{S}$ as well as the advection parameter $\gamma$. Given the simplified nature of the Stommel box model it is not possible to estimate appropriate values for these parameters from 
{first principles}. Data assimilation allows to reconstruct plausible values for these parameters from comparisons between model output and observations. In this section we will introduce the data assimilation methodology and the assimilated observations. The described model and methodology are implemented in the DA Python package DAPPER \cite{raanes_dapper_2023}.

\subsection{Ensemble Kalman filter}\label{data_assim}

{The true state of the system, which in our case consists of the temperature and salinity in the ocean boxes as well as the advection and diffusion coefficients, is not known. However, based on e.g. model output, it is possible to provide an estimate of the true value in the form of a probability distribution. Data assimilation (DA) is a procedure that combines this a priori distribution with imperfect observations (see \secRef{\ref{sec:obsproc}}) to produce a{n a} posteriori distribution with smaller spread around the truth. The most-likely value for the truth based on the a{ }priori distribution is referred to as the forecast, while the most-likely value of the {a} posteriori distribution is known as the analysis. In this work, we use a particular DA method known as the augmented ensemble transform Kalman (ETKF) \cite{bishop_adaptive_2001-1} in which both errors in the observations as well as the 
{a priori} distribution are assumed to be Gaussian. The mean and covariance of the latter estimated from an ensemble of model runs. A more detailed description of the augumented ETKF approach used can be found in appendix~\ref{app:kalman}.}

\subsection{Observations and Model Parameters}\label{sec:obsproc}

Geometry of the two boxes in the model, box-average temperature and salinity, observation error variances and box-surface-averaged temperature and salinity necessary for the model's forcing are compiled based on the Met Office's EN4 quality controlled objective analysis \cite{good_en4_2013} with bathymetric bias corrections \cite{gouretski_depth_2010}. The EN4 dataset is a dataset of several types of in-situ {potential temperature} and salinity measurements interpolated onto a 1x1 degree grid of the globe at 42 non-equally space depths on a monthly basis. From this dataset the vertical profiles in the North-Atlantic between {$23.5 \, ^\circ \mathrm{N}$} and $89.0 \, ^\circ \mathrm{N}$ and between $90.0 \, ^\circ \mathrm{W}$ and $90.0 \, ^\circ \mathrm{E}$ that are at least 1.5 km deep are selected. {These boundaries were chosen such model includes the North-Atlantic north of the intertropical convergence zone ($5-15 \,\unit{^\circ N}$) as well as Greenland and Norwegian . The southern boundary was chosen because precipitation in the former leads to freshening of the ocean surface that the simple Stommel box model geometry cannot accommodate. The northern boundary was chosen because overflow from the Nordic sea over the Iceland-Schotland ridge has been to be shown to be the dominant source for North-Atlantic Deep water (NADW) \cite{boning_deep-water_1996} which makes up the lower branch of the AMOC circulation. Only the part of the profiles above a depth of 3.5 km is used as this corresponds to the lower limit of occurrence of the NADW. Data from the period 2004-2022 are used. This selection is based on the fact that during this period extensive subsurface measurements were available thanks to the Argo program \cite{Argo} as well as that they exclude a sharp spike of surface temperatures at the beginning of the century probably caused by processes our simple Stommel box model cannot represent. 

The processes behind formation of North-Atlantic Deep Water (NADW) in the Artic that makes up the deep-ocean part of the AMOC are numerous, complex and their relative importance is disputed. They include the overflow of intermediate waters formed during the winter from the Nordic and Irminger seas over the Greenland-Scotland ridge \cite{boning_deep-water_1996,petit_atlantic_2020},  downward convection in the Labrador Sea driven by heat-loss \cite{bailey_formation_2005} and compensation for wind-driven upward mixing in the Atlantic \cite{kuhlbrodt_driving_2007,lozier_deconstructing_2010}. Of these only the former mechanism is incorporated in the Stommel model. Therefore we divide the vertical profiles between the boxes as follows: those profiles in which surface temperature drops below the depth-averaged temperature for on average 1 month per year or more, and consequently convection due to heatloss can take place, are assigned to the polar box. The other profiles are assigned to the equator box. This choice was not only motivated by the physical formation process of NADW, but also on more practical grounds. The data assimilation system needs to be able to reduce the average temperature of the deep ocean if it is found to exceed observed values. This definition of the polar box ensures a mechanism to release heat from the ocean is present. This choice is not entirely unproblematic (see section \ref{sec:concl}) and deviates from the one used by \cite{drijfhout_stability_2011,wood_observable_2019,chapman_quantifying_2024}. Their boxes extend to $23^\circ S$, cover the whole Artic ocean, not only the Nordic seas. Furthermore, they simple define the equatorial box to be the part of the Atlantic Ocean between $30^\circ S$ and $40^\circ N$ as they base deep-ocean temperatures and salinities on fit a to long-term climate model output. Hence overestimation of deep-ocean temperature can be dealt with by lowering initial deep-ocean temperatures and no heat release mechanism is necessary. Furthermore, the fit in \cite{chapman_quantifying_2024} is based solely on temperature and salinity differences between the boxes and hence the box model does not need to provide a deep ocean heat release mechanism. In other box model studies \cite{monahan_stabilization_2002} no attempt is even made to match boxes and their values with actual observations.}

\begin{table}[h!]
	\centering
	\begin{tabular}{|c|cc|}
		\hline
		& Polar & Equator \\
		\hline
		$\dx$ $\mathrm{[km]}$ & 7274 & 7274 \\
		$\dy$ $\mathrm{[km]}$ & 324 & 2262 \\
		$\dz$ $\mathrm{[m]}$ & 3148 & 3148\\
		\hline
	\end{tabular}
	\caption{\raggedright The zonal width, meridional width and depth of the two boxes used in the Stommel model based the EN4 data. }
	\label{tab:geometry}
\end{table}

\begin{figure}
	\centering
	\includegraphics[width=.6\textwidth]{ 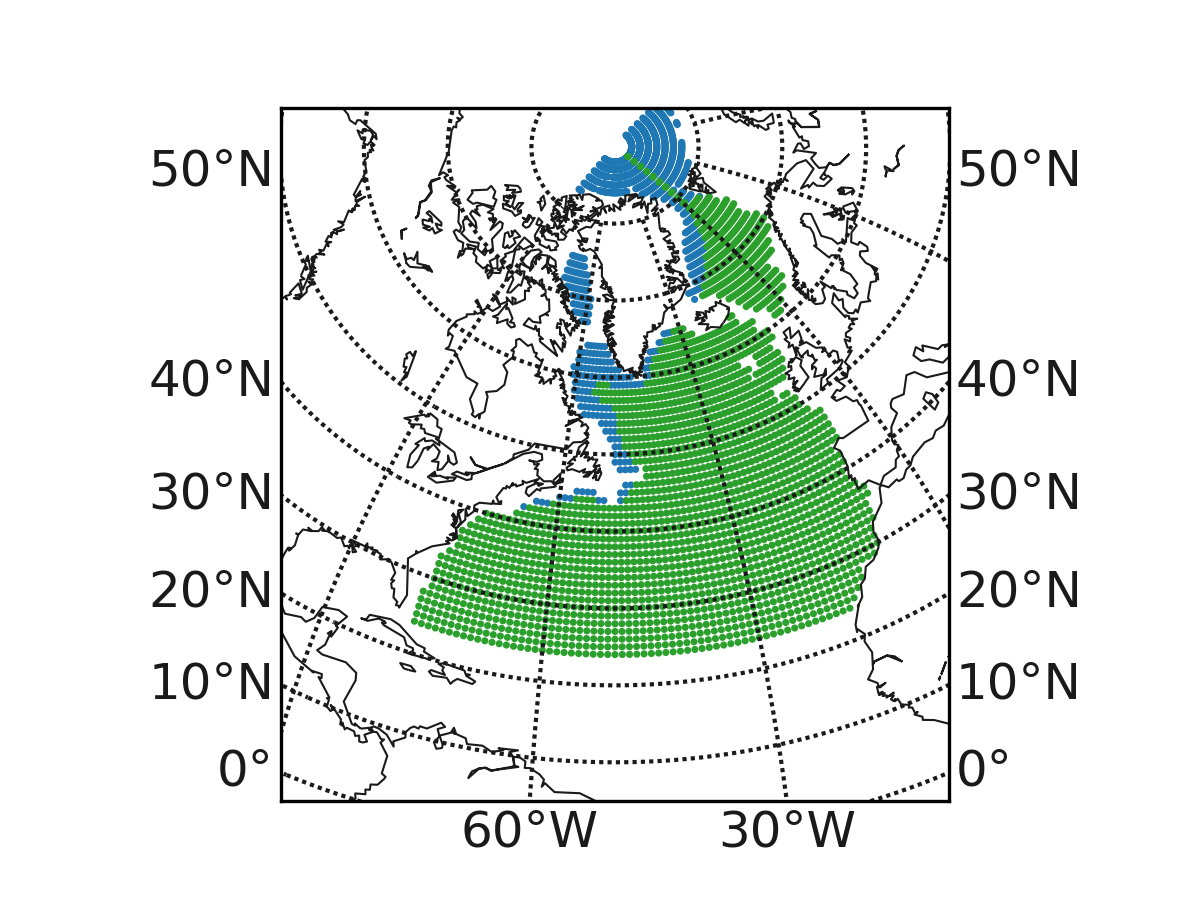}
	\caption{\raggedright An example of EN4 profile assignment in the year 2010. Green points signify profiles assigned to the equatorial box, and blue points signify profiles assigned to the polar box.}
	\label{fig:geometry}
\end{figure}

Volume $V$ of each of the boxes is calculated by summing the volume of all grid cells in the box. Similarly, surface area $A$ is calculated by summing the area of the EN4 grid cells in the top layer in both boxes. The depth of both boxes is then calculated as $\Delta z = \frac{V_{p}+V_{e}}{A_{p}+A_{e}}$. The cross-section $A_{c}$ between boxes is derived by summing the vertical area between profiles assigned to the polar box and those assigned to the equatorial box. An example assignment is shown in Figure \ref{fig:geometry}. From this zonal and meridional size of the boxes is defined as $\Delta x = \frac{A_{c}}{\Delta z}$ and $\Delta y_{p}=\frac{A_{p}}{\Delta x}$, $\Delta y_{e}=\frac{A_{e}}{\Delta x}$. The values resulting from this exercise are shown in \tabRef{\ref{tab:geometry}}. 

\begin{table}[h!]
	\centering
	\begin{tabular}{|c|cc|}
		\hline
		& Time-averaged obs. & Obs. error std. dev. \\
		\hline 
		$T_{p*}$ $[\celsius]$ & 1.2 &  0.3 \\
		$T_{e*}$ $[\celsius]$ & 5.5 & 0.5 \\
		$S_{p*}$ $[\ppt]$ &  34.83 & 0.07 \\
		$S_{e*}$ $[\ppt]$ & 35.15 & 0.07 \\
		\hline
	\end{tabular}
	\caption{\raggedright Time-average of the assimilated temperature and salinity observations together with the observational standard deviation, i.e. the square-roots of the diagonal of $\Sigma_{d}$. }
	\label{tab:TS_box}
\end{table}

For each month $t$ and each of the two boxes the box-averaged temperature is calculated as 
\begin{equation}
	T_{\cdot *}(t) = \big(\sum_{i} V_{i} T_{i,*}(t)\big)\big(\sum_{i} V_{i}\big)^{-1},
	\label{eq:T_mean}
\end{equation}
where $i$ runs over all EN4 grid cells in the box below the surface layer, $V_{i}$ is the volume of the grid cell, $\sigma_{T,i}$ is the uncertainty in temperature as specified in the EN4 data set. The variance in the value is calculated as 
\begin{equation}
	\sigma_{T*}^{2}(t) = \big(\sum_{i} V_{i} \sigma^{2}_{T_{i,*}}(t)\big)\big(\sum_{i} V_{i}\big)^{-1}.
	\label{eq:T_sigma}
\end{equation}
These calculations are repeated for salinity $S$.  The values after January 2004 are used to create the observations to be assimilated in $\Vec{d}(t)=[T_{p*}(t),T_{e*}(t),S_{p*}(t),S_{e*}(t)]^{\mathrm{T}}$. The $4 \times 4$ observational error covariance matrix $\Sigma_{d}$ is assumed to be diagonal. The values on the diagonal are chosen to be the maximal values for $\sigma_{T*}^{2}(t)$ and $\sigma_{S*}^{2}(t)$ in each box. The resulting time-average values for the observations as well as the observational error standard deviation found this way are shown in \tabRef{\ref{tab:TS_box}}.

\begin{table}[h!]
	\centering
	\begin{tabular}{|c|cccc|}
		\hline
		Field & $\cdot_{0}$ & $\cdot_{\cos}$ & $\cdot_{\sin}$ & seasonal amplitude \\
		\hline 
		$T^{a}_{p}$ $[\celsius]$ & 1.5 & -1.5 & -1.1 & 1.9\\
		$T^{a}_{e}$ $[\celsius]$ & 16.7 & -2.4 & -2.3 & 3.4 \\
		$S^{a}_{p}$ $[\ppt]$ & 33.05 & 0.22 & 0.32 & 0.39 \\
		$S^{a}_{e}$ $[\ppt]$ & 35.77 & 0.04 & 0.05 & 0.07 \\
		\hline
	\end{tabular}
	\caption{\raggedright Regression coefficients found by fitting model \eq{\ref{eq:T_regression}} to the surface observations for temperature $^{a}$ and salinity $^{a}$ in polar $\cdot_{p}$ and $\cdot_{e}$ boxes. Also shown is the seasonal amplitude $\sqrt{\cdot_{\cos}^2+\cdot_{\sin}^2}$ }
	\label{tab:TS_surface}
\end{table}

The calculations in \eq{\ref{eq:T_mean}} and \eq{\ref{eq:T_sigma}} are repeated but now with $i$ running over all grid cells in the surface layer of a box. The signal is then represented as 
\begin{equation}
	T^{a}(t) \approx T^{a}_{0} + T^{a}_{\sin} \sin(\frac{2\pi t}{\tau}) + T^{a}_{\cos} \cos(\frac{2\pi t}{\tau}) ,
	\label{eq:T_regression}
\end{equation}
with $\tau=1 \, \mathrm{year}$ and $T_{0}$, $T_{\cos}$, $T_{\sin}$ determined by a weighted linear regression against the box-averaged surface values using $\sigma_{T*}^{-2}(t)$ as weights. The regression coefficient for surface temperature and their equivalents for surface salinity are given in \tabRef{\ref{tab:TS_surface}}.

\begin{table}[h!]
	\centering
	\begin{tabular}{|c|ccc|}
		\hline 
		& dimensional & scale & nondimensional \\
		\hline 
		$\Omega$ & $2\pi$ $\mathrm{rad\,yr^{-1}}$ & $\frac{k_{T}}{\dz}$ & $3564$\\
		$B$ & 1.51 $\celsius$ & $\frac{k_{T}}{\alpha_{T} \gamma \dz}\frac{\dy_{p} \dy_{e}}{\dy_{p} + \Delta y_{e}}$ &
		$1.18$ \\
		$\hat{B}$ & 0.32 $\unit{\ppt}$ & $\frac{k_{T}}{k_{S}}\frac{k_{T}}{\alpha_{S} \gamma \dz} \frac{\dy_{p} \dy_{e}}{\dy_{p} + \dy_{e}}$ & $0.66$ \\
		$A$ & - & - & $0.52$\\
		\hline
	\end{tabular}
	\caption{\raggedright Relation between the coefficients $\Omega$, $A$, $B$, $\hat{B}$ in \cite{budd_dynamic_2022} and the surface forcing obtained from the EN4 observation data.}
	\label{tab:buddParam}
\end{table}

The regression coefficients in \tabRef{\ref{tab:TS_surface}} can be related to the perturbations added to $\eta_{1}$ and $\eta_{2}$ in Equation~\eqref{eq:eta1_B} and Equation~\eqref{eq:eta2_Bhat}. These relations can be found in \tabRef{\ref{tab:buddParam}}. 

\begin{table*}[htbp]
	\caption{Parameters, forcings, initial conditions and their sources. 
	}
	\resizebox{\textwidth}{!}{%
		\begin{tabular}{|l|l|l|l|}
			\hline
			\bf{Data} & \bf{Values} & \bf{Source} & \bf{Details} \\ \hline
			Temperature (ocean temperature)    & Pole: $T^{*}_{p}=1.23 \,\unit{\celsius}$, Equator: $T^{*}_{e}=5.4 \,\unit{\celsius}$ 
			&  &  Initial values from EN4 data
			\\ \hline
			Salt (ocean salinity)    & Pole: $S^{*}_{p}=34.82 \,\unit{ppt}$
			Equator: $S^{*}_{e}=35.15 \,\unit{ppt}$
			&  &  Initial values from EN4 data
			\\ \hline
			temp\_diff (surface heat mixing coefficient)    & $k_{T}=3.7\cdot 10^{-6} \,\unit{m s^{-1}}$
			&  & equation \ref{eq:minimization}
			\\ \hline
			salt\_diff (surface salinity mixing coefficient)    & $k_{S}=1.2\cdot 10^{-6} \,\unit{m s^{-1}}$
			&  & equation \ref{eq:minimization}
			\\ \hline
			$\rho_0$ (reference density)   & $1027 \,\unit{kgm^{-3}}$
			& \cite{Nyc} & 
			\\ \hline
			$S_0$ (reference salinity)   & $35 \,\unit{ppt}$
			& \cite{Nyc} & 
			\\ \hline
			$T_0$ (reference temperature) & $10 \,\unit{\celsius}$
			& \cite{Nyc} & 
			\\ \hline
			$\alpha_T$ (thermal expansion coefficient) & $\frac{0.15 \,\unit{\celsius^{-1}}}{\rho_{
					ref}}$
			& \cite{Nyc} & 
			\\ \hline
			$\alpha_S$ (haline expansion coefficient) & $\frac{0.78 \,\unit{ppt^{-1}}}{\rho_{
					ref}}$
			& \cite{Nyc} & 
			\\ \hline
			Initial estimate $Q_{overturning}$ (meridional overturning flux)) & $18 \,\unit{Sv}$
			& \cite{MCC} & 
			\\ \hline
			$\gamma$ (advective transport coefficient) & $2.0  \,\unit{ms^{-1}}$
			&  & equation \ref{eq:minimization}
			\\ \hline
			$T^a(t)$ & Oscillates over time
			&  &Calculated from EN4 dataset;  see table \ref{tab:TS_surface}
			\\ \hline
			$S^a(t)$ & Oscillates over time
			&  &  Calculated from EN4 dataset; see table \ref{tab:TS_surface}
			\\ \hline
			V\_ice & $2.9 \times 10^6 \,\unit{km^3}$
			&  & 
			\\ \hline
			Greenland ice sheet melt period $t_{melt}$ & $10000 \,\unit{y}$ & \cite{khan_greenland_2015,mouginot_forty-six_2019,mankoff_greenland_2020}&
			\\ \hline
			Equatorial warming rate & $0.03 \,\unit{\celsius yr^{-1}}$  & \cite{clim, clim2, you_warming_2021} &
			\\ \hline
			Polar warming rate
			& $0.06 \,\unit{\celsius yr^{-1}}$  & \cite{Zha,you_warming_2021} & 
			\\ \hline
			Model time step & 1 month & & 
			\\ \hline
		\end{tabular}
		\label{tab:param}
	}
\end{table*}

The numerical scheme for Stommel model uses a finite volume approach with an explicit 4th-order Runge-Kutta scheme with a time step of 1 month. Initial conditions for the ocean temperatures and salinities in the different ensemble members are drawn from a Gaussian distribution having the EN4 box-averaged temperatures and salinities for January 2004 as mean and $\Sigma_d$ as covariance. {Initial values for the logarithm of the diffusion and advection parameters are estimated by minimizing, using Scipy's Nelder-Meat method,
\begin{equation}
 \frac{(\Delta T_{*}-\Delta T_{eq})^{2}}{\sigma_{T,e}^2+\sigma_{T,p}^2}+\frac{(\Delta S_{*}-\Delta S_{eq})^{2}}{\sigma_{S,e}^2+\sigma_{S,p}^2}+\frac{(Q_{*}-Q_{eq})^{2}}{(2.5\,\unit{Sv})^2}
 \label{eq:minimization}
\end{equation}
with $\Delta T_{*}$, $\Delta S_{*}$ the difference between equatorial and polar box's initial ocean temperature and salinity, $\sigma_{T,\cdot}$, $\sigma_{S,\cdot}$ the observational standard deviations in the 3rd column of table~\ref{tab:TS_box}, $Q=18 \,\unit{Sv}$ the observed meridional transport and $\Delta T_{eq}$, $\Delta S_{eq}$ $Q_{eq}$ equilibrium values for temperature, salinity difference and transport that depend on the model parameters. I.e. we pick parameters such that our initial conditios are close to an equilibrium value.} 
The logarithms are perturbed by adding realizations from a Gaussian with a standard deviation of $0.26$. This corresponds to the assumption that the relative error in these model parameters is $30\%$. {These model parameters are assumed to be constant in time. More specifically, their values are held constant between DA corrections during the prediction step (see appendix~\ref{app:kalman}. However, the parameters for each ensemble member might still change over time due to the DA corrections reflecting that the probability density distribution for the true value of the model parameters is not necessarily stationary.} Values for the other parameters used in the model can be found in Table~\ref{tab:param}. {Note that this table includes an equatorial and a polar warming rate, and that the polar rate is twice the equatorial. The choice to use two warming rates was made to account for Arctic amplification. Arctic amplification is the name for an effect observed over many decades in which the Arctic warmed faster than other parts of the globe\cite{walsh2014intensified}. The precise magnitude of the difference is disagreed upon, but many studies report the poles warm around twice the global average \cite{richter2020state,jansen2020past}, so we used this ratio to inform our model.}

\section{Twin Experiments}\label{sec:twin}

In this section we will evaluate whether the data assimilation system is capable of correcting ocean temperatures and salinities as well as reconstructing the Stommel model parameters. In addition to this we want to see what happens with the AMOC circulation if the effects of climate change are added to the model. {This is done by running an Observing System Simulation Experiment (OSSE), or twin experiment \cite{halliwell_rigorous_2014}. In such an experiment one of the ensemble members is set apart as ``truth'' and withheld from the DA system forming one of the twins. Artificial observations are created from the ``truth'' by adding perturbations drawn from a Gaussian distribution with covariance $\Sigma_{d}$ with these perturbations mimicking observation errors. After assimilating these observations into the experiment, the forecasts and analyses, making up the other twin, can be compared with the ``truth'' to quantify performance of the DA system.} 

\begin{figure*}[tbp]
	\centering
     \includegraphics[width=.7\linewidth]{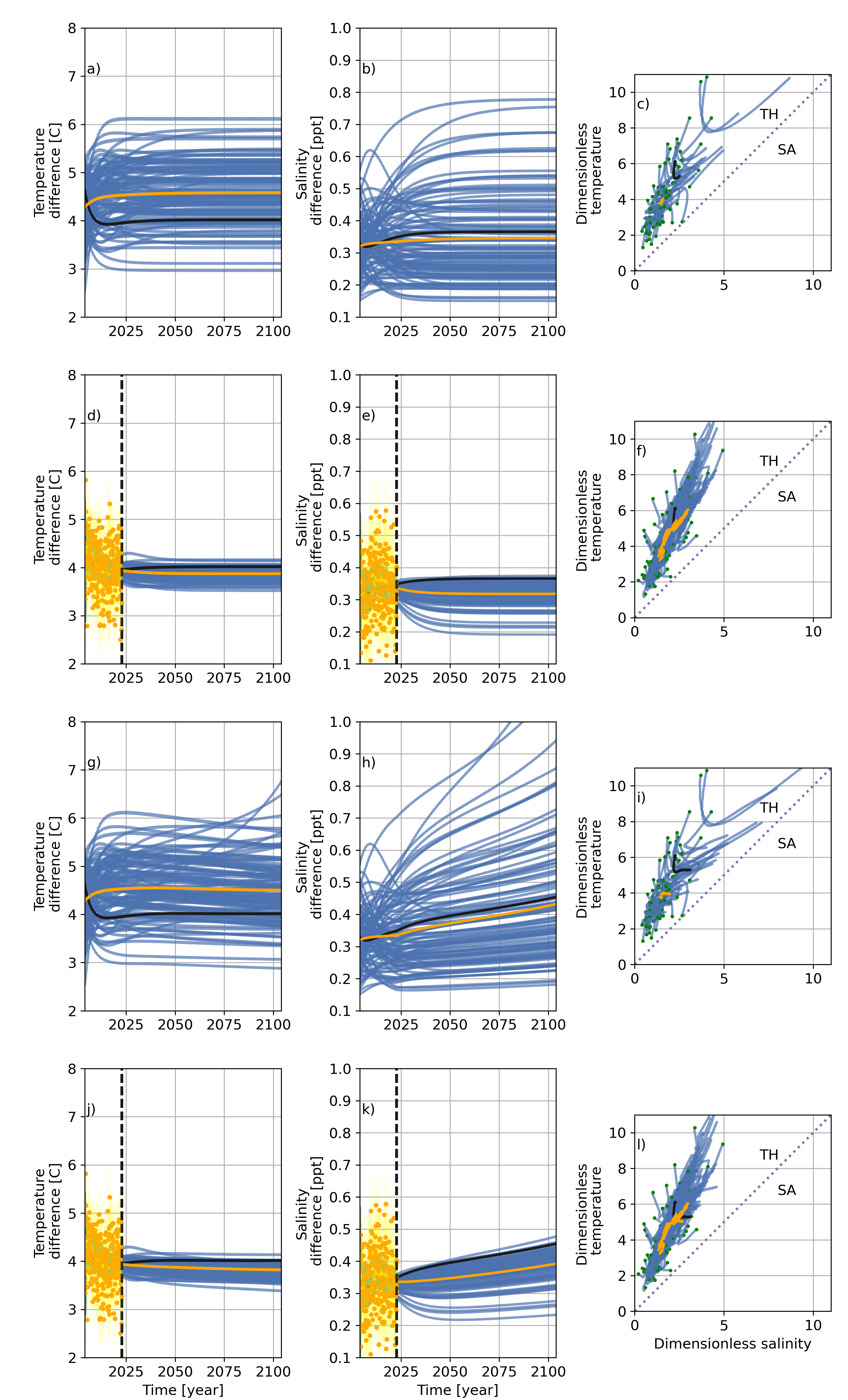}
    \caption{\raggedright Forecast time series of temperature and salinity differences between equatorial and pole boxes for 100 ensemble members, with corresponding dimensionless salinity-vs-temperature phase portraits. Green dots: initial conditions; dashed blue line: TH-SA regime boundary; black paths: synthetic truth; orange paths: ensemble mode. Rows show different experiments: a)-c) nCC-nDA, d)-f) nCC-yDA, g)-i) yCC-nDA, j)-l) yCC-yDA. (nCC: no climate change, yCC: with climate change, nDA: no data assimilation, yDA: with data assimilation). Columns: temperature difference, salinity difference, dimensionless temperature vs. salinity phase portrait.}
	\label{fig:TSTwin}
\end{figure*}

We conduct simulations under four distinct scenarios. In the first scenario initial conditions and model parameters are perturbed in the different ensemble members but no data assimilation is used. The second scenario mirrors the first but now ocean temperatures and salinities are assimilated each month during the period 2004-2022. For the third scenario the effects of global warming are taken into account. In this study we take two processes related to climate change into account. First, that the average temperature increases faster near the pole  ($0.06 \unit{\frac{\celsius}{year}}$ for polar box) than at lower latitudes ($0.03 \unit{\frac{\celsius}{year}}$ for equatorial box). Second, that the Greenland ice sheets melts add a constant rate $m^{*}=\frac{V_{ice}}{t_{ice}}$ with $V_{ice}$ and $t_{ice}$ given in Table~\ref{tab:param}.  In order not to interfere with the DA, warming and melting only start after 2022, i.e. after the end of the DA period, and is applied uniformly to all ensemble members. The fourth scenario matches the third and includes data assimilation. The model time step is 1 month and DA is carried out monthly. \figRef{\ref{fig:TSTwin}} shows the model ocean temperatures and salinities over 100 years for 100 ensemble members together with its most-likely value (the ensemble mean) and the artificial truth. {For times at which an {\it update step} is carried out (see appendix~\ref{app:kalman}) and two solutions are available, the solution after the DA correction is shown. The timeseries obtained in this way is referred to as the forecast solution.}
Temperatures and salinities in the 3rd column are non-dimensionalized using the scaling constants in Section~\ref{sec:stommel_model}. Blue curves indicate the thermohaline (TH) regime, while red curves represent the salinity-driven (SA) regime. {As shown by \eqref{eq:dimensionless}, the solution to the Stommel equations can be written solely as function of the difference in temperature and salinity between the equatorial and polar box. Therefore the results in \figRef{\ref{fig:TSTwin}} are depicted as function of these differences only.} 

As seen in Figure \ref{fig:TSTwin}, none of the ensemble members predict a change in flow direction in the experiment with no climate change and no data assimilation. This is notwithstanding the large uncertainties in the model parameters. One can see that there is a significant spread between trajectories in the time plots, but the ensemble remains monomodal, i.e. the ensemble members do not form different clusters. This is a necessary condition for the EnKF to function properly. Comparison of \figRef{\ref{fig:TSTwin}a,b,g,h} with \figRef{\ref{fig:TSTwin}d,e,j,k} respectively shows that assimilating temperatures and salinities for the period 2004-2022 reduces the spread bringing all ensemble members closer to the truth. Trajectories become and stay tightly connected in the time and phase plots. This indicates that DA is successful in correcting the temperatures and salinities. 
When the effects of climate change are applied in Figures \ref{fig:TSTwin}g,h,j,k the solutions no longer converges towards an equilibrium value. Instead temperature and salinity keep increasing reflecting the introduced changes in forcing. Apart from this, the results are qualitatively the same as in the first two {scenarios}: the circulation remains in TH mode and the DA reduces the errors in the ensemble members.  

\begin{figure*}[tbp]
	\centering
    \includegraphics[width=\linewidth]{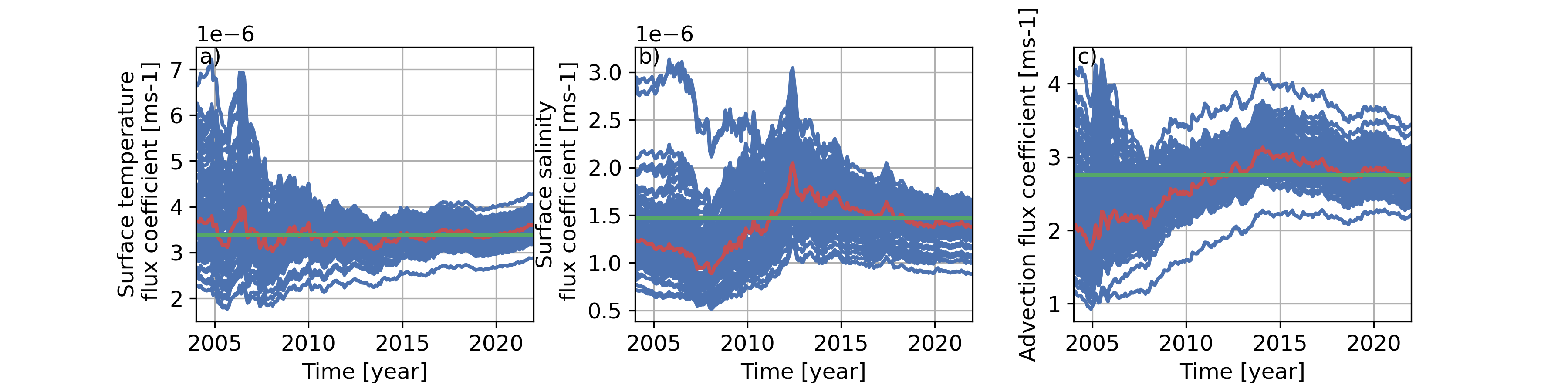}
	\caption{\raggedright  Read left to right: a) temperature ($\kappa_{T}$), b) salinity ($\kappa_{S}$), and c) advective flux coefficients ($\gamma$) of the ensemble plotted over the assimilation period for an experiment assimilating 18 years of synthetic data without any global warming effects (the experiment with global warming is effectively identical). The red curve represents the ensemble mode, while the green line represents the true coefficient value.}
	\label{fig:advectiveTE}
\end{figure*}

In order to verify whether the DA system is capable of correcting the model parameters, temperature diffusivity parameter, salinity diffusivity parameter and advection parameter are shown in Figure~\ref{fig:advectiveTE} as function of time during the DA period together with their true value. We see that the spread is sound as so far as the truth lies within the range of the ensemble. Comparing the most-likely estimate with the truth (in green), we see that initially DA can increase the error in the model parameters. However, by the end of the period errors have reduced to values below the initial error. 

\begin{figure*}
	\centering
    \includegraphics[width=.75\linewidth]{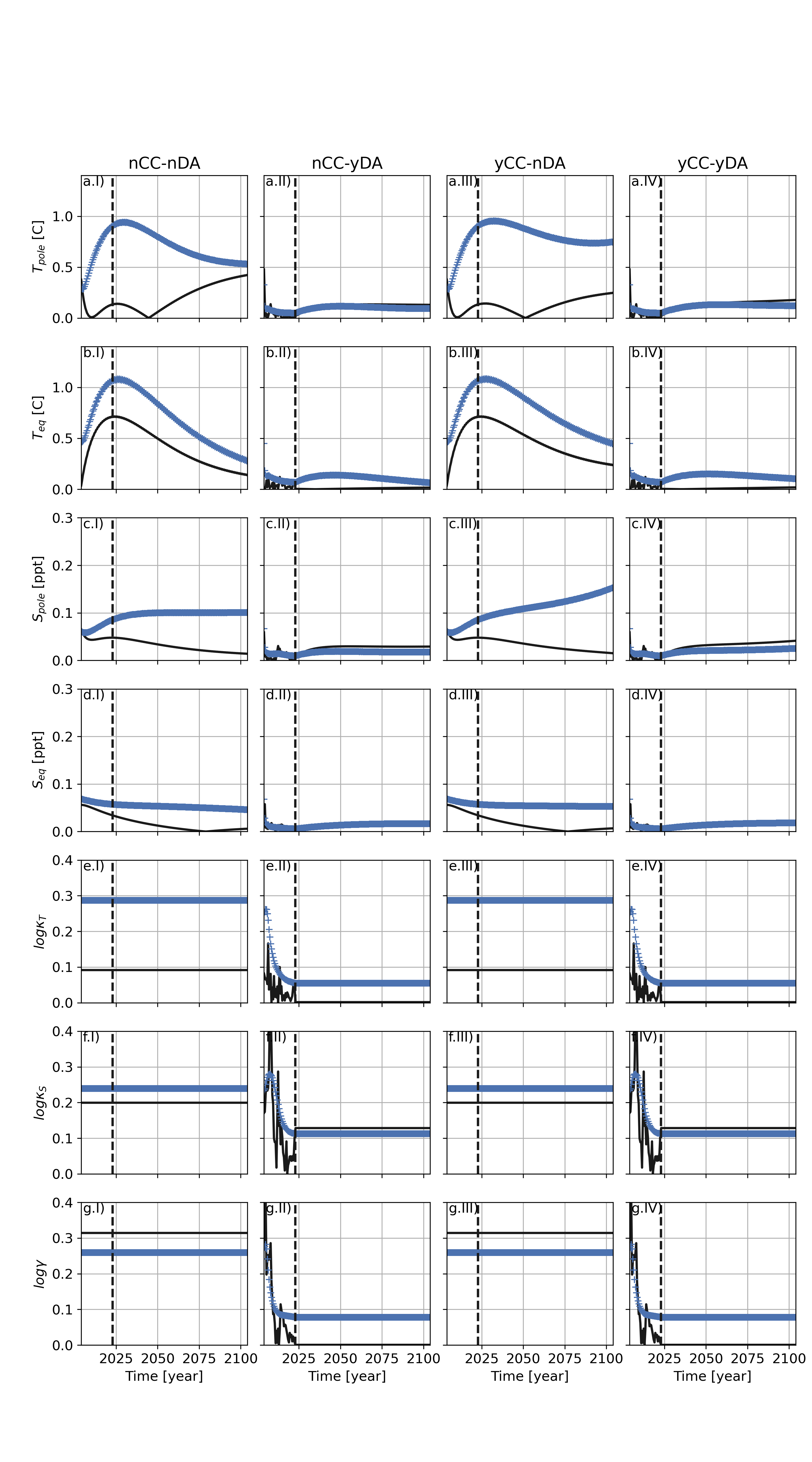}
	\caption{\raggedright Read left to right, top to bottom. (Black solid lines) RMSE and (blue pluses) ensemble standard deviation of (rows) different variables as function of time for (columns) the four different experiments in Figure \protect\ref{fig:TSTwin}. For visibility only every 6th month is shown.}
	\label{fig:RelTwin}
\end{figure*}

To quantify the behavior of the ensembles in \figRef{\ref{fig:TSTwin}} in more detail, the ensemble standard deviation, or ensemble spread, and RMSE for the different experiments are shown in \figRef{\ref{fig:RelTwin}} for the different state variables. The first thing that stands out is that ensemble spread and RMSE after the DA period generally stay constant or decrease, even in the absence of DA. The second noticeable point is that DA is {successful} in reducing the error in the model parameters. Errors for the advection parameter (\figRef{\ref{fig:RelTwin}g}) and temperature diffusion parameter (\figRef{\ref{fig:RelTwin}e}) are decimated while that of the salinity parameter is reduced but not to the same extent (\figRef{\ref{fig:RelTwin}f}). This reflects that salinity diffusion parameters is smaller. Consequently, it takes the system more time to equilibriate with surface salinity forcing and consequently the DA is less sensitive to changes in the salinity diffusion coefficient than in the advection and temperature diffusion coefficient. The fact that by the end of the DA period the spread in the salinity diffusion parameter is still decreasing (\figRef{\ref{fig:RelTwin}l}), whilst that in the temperature diffusion parameter (\figRef{\ref{fig:RelTwin}j}) and advection parameter ((\figRef{\ref{fig:RelTwin}n}) has leveled off suggest that further RMSE reduction for the former could be achieved by extending the DA period. Finally, for a properly calibrated ensemble Kalman filter the expectation value of the ensemble spread and root mean-square error (RMSE) should be equal. Comparison of the RMSEs and ensemble standard deviations at the end of the DA period shows they are of similar magnitude with the spread slightly overestimating the uncertainty. Thus confirming that DA is performing as expected.

\begin{figure*}[tbp]
	\centering
	\includegraphics[width=.8\linewidth]{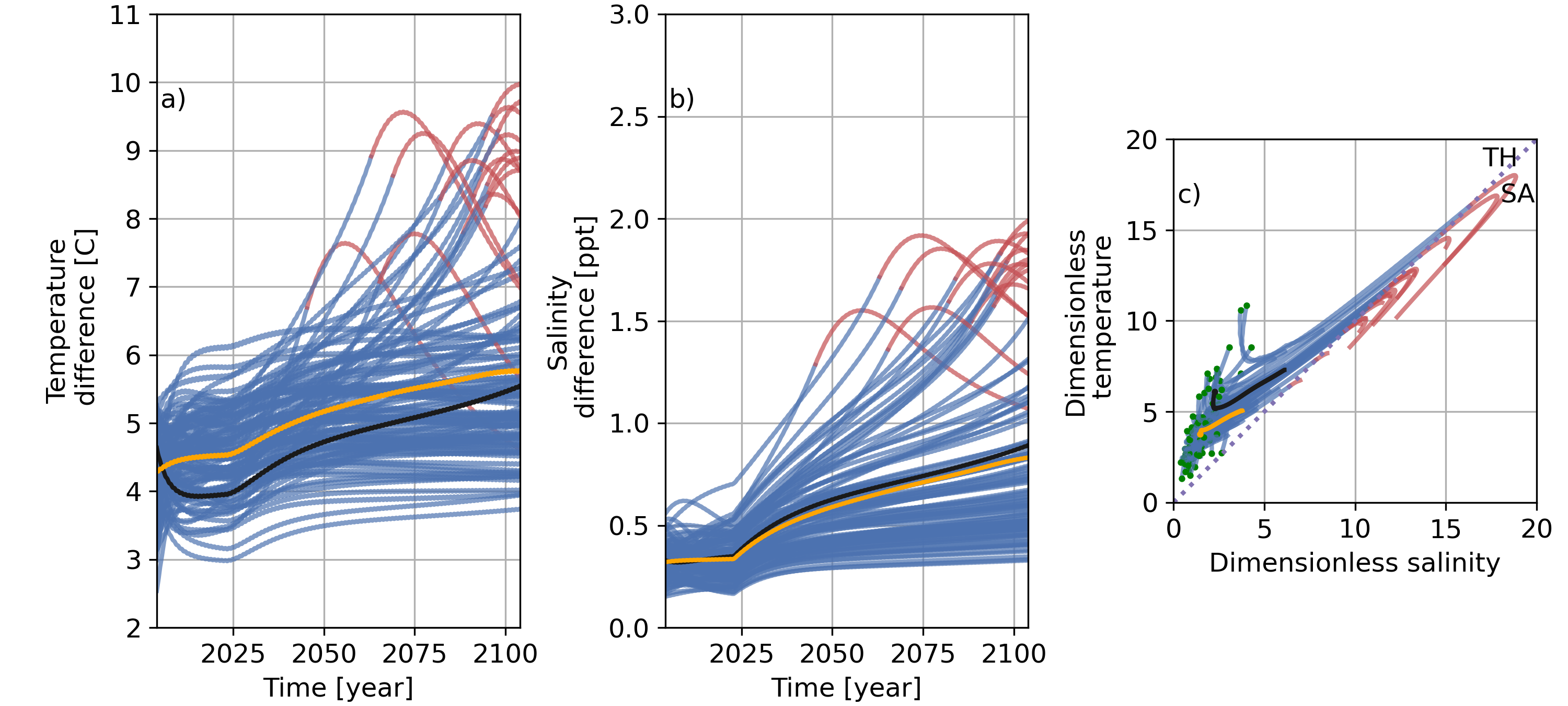}
	\caption{\raggedright  Read left to right:  a) temperature and b) salinity time series, as well as c) a phase portrait for an experiment with yearly warming of $0.03 \,\unit{\celsius}$ per year in the equator ($0.06 \,\unit{\celsius}$ per year at the pole) and a melt period of 1000 years without DA. Trajectories turning red indicate a predicted flipping of flow direction. The black trajectories represent the synthetic truth values, and the orange represents the ensemble mode.}
	\label{fig:twin_flip}
\end{figure*}

All ensemble members in \figRef{\ref{fig:TSTwin}} are in the TH circulation mode. In order to see whether the setup is capable of producing solutions with the SA circulation, the impacts of the climate change are scaled up. More precisely, the warming rate in the equatorial box is increased to $0.07 \,\unit{\celsius}$ per year (the polar box warming rate is raised to $0.14 \,\unit{\celsius}$ per year) and the melting period for the Greenland ice sheet is reduced from 10,000 to 1,000 years. Once again, these climate change driven effects are only switched on in 2022  otherwise the DA system would try to correct errors in the forcing by modifying the model parameters. In \figRef{\ref{fig:twin_flip}b} salinity  between-box difference in all ensemble members increases due to the effects of climate change. This reduces the circulation and temperature difference start to increase as consequence of the temperature difference present in the surface. Eventually, \figRef{\ref{fig:twin_flip}a,b} show that depending on initial conditions and the model parameters of the ensemble member, there is a rapid increase in the temperature and salinity between-box difference resulting in a reversal of the AMOC circulation direction. So, the setup is capable of producing solutions with a circulation change.  

To summarize, in this section we tested the ability of our setup to predict the truth using a twin experiment. We find that the setup produces an ensemble that covers the truth and can, with the sufficient forcing, produce solutions both with the TH as well as the SA circulation. DA is successful in reducing the spread and errors in both temperatures, salinities, diffusivity parameters and advection parameter. The 18 year DA period used is sufficient to reduce the uncertainty in temperature diffusivity and advection parameter to the lowest value attainable with this setup, but the uncertainty in the salinity diffusivity could have been reduced further by extending the DA period. Having established that our DA setup has the ability to recover the model parameters from temperature and salinity observations, we will use those observations to find realistic values for these parameters in Section~\ref{sec:obs}. 

\section{Assimilating Box-Averaged observations in the Non-Autonomous Stommel Model}\label{sec:obs} 

In this section, we run simulations using \eqref{eq:nonaut} while assimilating box-average observations compiled from the EN4 dataset as described in Section~\ref{sec:obsproc}. In this section we pursue a two-fold aim. First, we attempt to find realistic values for the salinity diffusivity $k_S$, temperature diffusivity $k_T$ and advection parameter $\gamma$ as well as the non-dimensional values $\eta_{1}$, $\eta_{2}$ and $\eta_{3}$ which depend on these parameters. Second, having found realistic values, we continue running the box model forward in time to determine the {likelihood} of a circulation reversal happening in this century both for the case with climate change and the case without.

\subsection{Interfered Model Parameters}
\label{sec:nonauto}

\begin{figure*}
	\center
	\includegraphics[width=\linewidth]{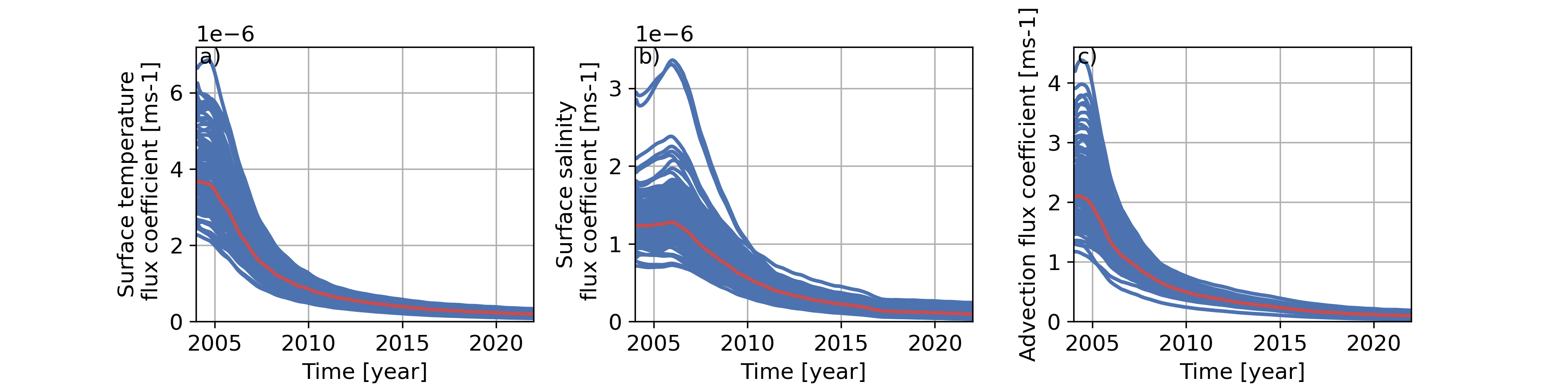}
	\caption{\raggedright Read left to right: a) surface temperature ($\kappa_{T}$), b) surface salinity diffusivity ($\kappa_{S}$) and c) the advection flux coefficients ($\gamma$) for the different ensemble members as functions of time.  The red curves represents the ensemble mode. \label{fig:realParams}}
	
\end{figure*}

\begin{table}[h!]
	\centering
	\begin{tabular}{|c|cc|}
		\hline 
		Parameter &  Most-likely & 90\%-confidence interval \\
		\hline 
		$k_{T}$  $[10^{-7} \mathrm{ms^{-1}}]$ &  $1.8$ &  $[1.2, 2.9]$\\
		$k_{S}$  $[10^{-7} \mathrm{ms^{-1}}]$ &  $0.8$  & $[0.5, 2.1]$ \\
		$\gamma$ $[10^{-2} \mathrm{ms^{-1}}]$ & $8.5$ & $[6.1, 13.9]$ \\
		\hline
	\end{tabular}
	\caption{\raggedright Most-likely estimates of model parameters after assimilating monthly boxed-averaged temperatures and salinities for the period 2004-2022 together with their $90$\%-confidence interval.}
	\label{tab:bestParams}
\end{table}

The three different model parameters are shown in \figRef{\ref{fig:realParams}} as function of time. As the effects of climate change are never applied during the DA period (see \secRef{\ref{sec:twin}}), the figure will look the same for both experiments with and without climate change. Notice that we are still assuming that within the conceptual constraints of the Stommel box model the true value of these parameters is fixed in time. The time dependence reflects that the probability distribution for the truth is changing as more and more observations get assimilated. Had we used the Rauch-Tung-Striebel smoother \cite{rauch_maximum_1965} we would have been able to correct the parameters in the beginning of the DA window using observations at the end of the window. However, as this would have not changed the final estimate, i.e. the value for the model parameters at the end of the DA window used to propagate the ensemble members forward, we have abstained from the additional effort. The figure shows an initial overestimation of the value of all three coefficients: the mode of the ensemble and all ensemble members exhibit a decreasing trend that starts to level off near the end of the DA period. Similar to the twin experiment in \secRef{\ref{sec:twin}} the salinity diffusion parameter is slower to converge to its final value. The ensembles for all three parameters show a steadily decreasing spread indicating that the uncertainty of the parameter's true values decreases. The final estimates for the parameters, i.e. their values at the end of the DA period, together with their 90\% confidence interval based on the ensemble spread are given in \tabRef{\ref{tab:bestParams}}.

\begin{figure}[tbp]    
	\centering
	\includegraphics[width=.5\textwidth]{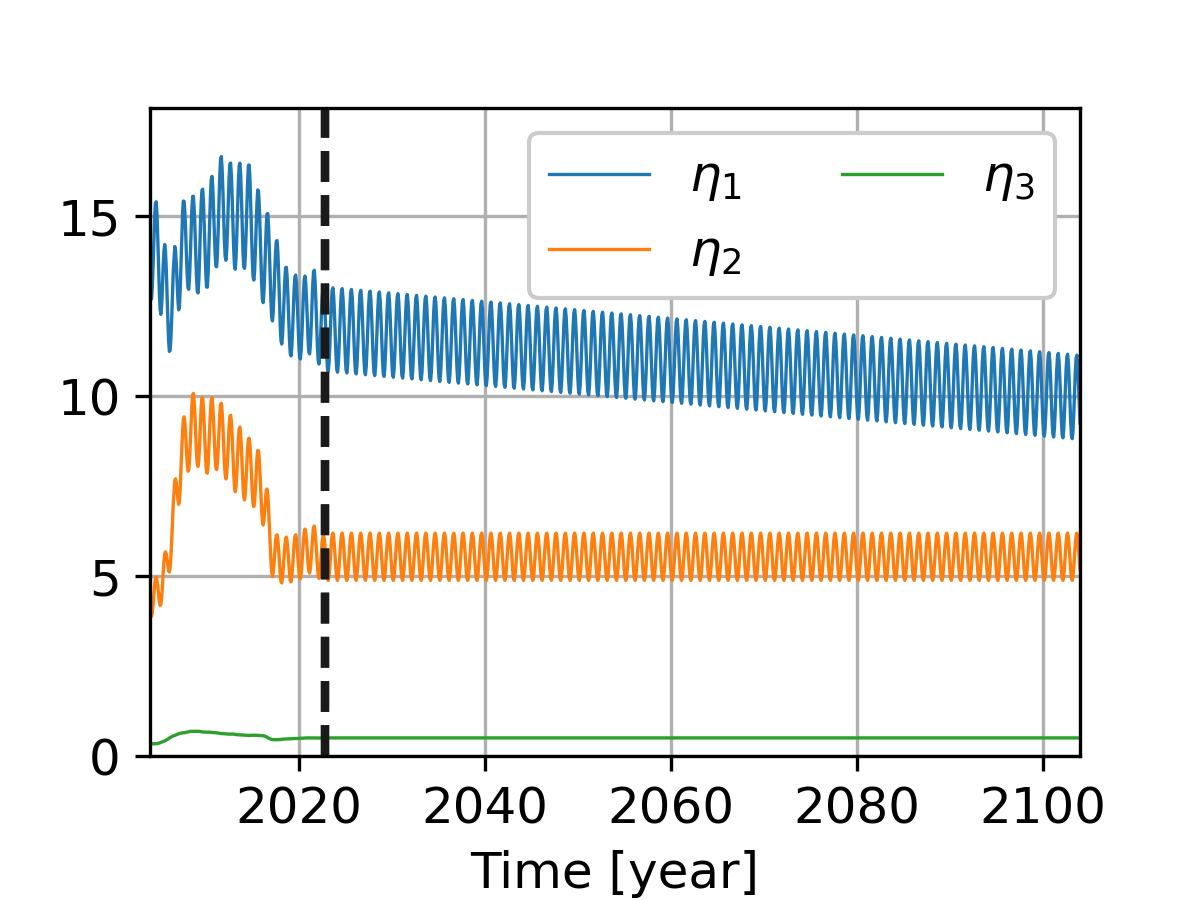}
	\caption{\raggedright The values of parameters $\eta_{1}$ (Equation~\eqref{eq:eta1_B}), $\eta_{2}$ (Equation~\eqref{eq:eta2_Bhat}) $\eta_{3}$ (Equation~\eqref{eq:eta3}) as function of time for the experiment with climate change. \label{fig:realEta}}
\end{figure}

The values for the non-dimensional parameters $\eta_{1}$, $\eta_{2}$ and $\eta_{3}$ as appearing in Equation~\ref{eq:dimensionless}
are shown in \figRef{\ref{fig:realEta}}. The annual oscillations clearly visible in $\eta_{1}$ and $\eta_{2}$ are the result of seasonal cycle imposed on the surface forcing as Equation~\eqref{eq:T_regression}. As climate change warms the polar box faster than the equatorial one, the temperature difference between the two decreases over time resulting in a decreasing slope in $\eta_{1}$ after the DA period. The annually averaged values for the ensemble mode of $\eta_{1}$, $\eta_{2}$ and $\eta_{3}$ at the end of the DA period are $9.2$, $4.2$ and $0.5$ respectively. These are about double the values explored in other Stommel box model literature \cite{dj,budd_dynamic_2022}. The absence of certain physical relevant processes for the AMOC and might explain difference. See also the discussion of transport at the end of this section. 
\section{Future Circulation}\label{sec:future}

\begin{figure*}[tbp]
	\centering
	\includegraphics[width=.8\linewidth]{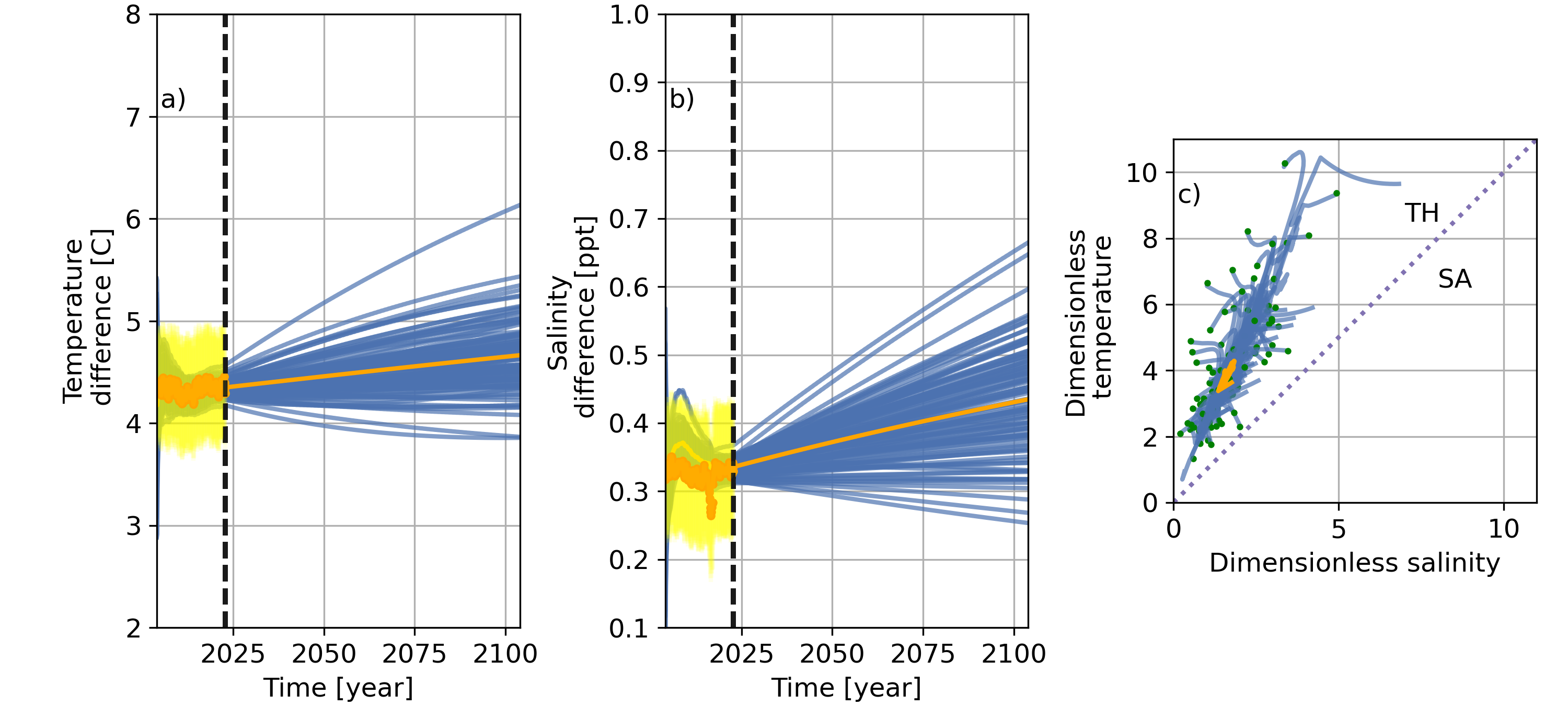}
	\caption{\raggedright Read left to right. Forecast time series starting in 2004 representing the equator-pole difference in a) temperature and b) salinity for the different ensemble members as well as c) the ensemble trajectories in a dimensionless temperature vs salinity phase portrait. The orange dots with bright yellow errorbars represent the observed values and their standard deviation. Solutions with a TH circulation are shown in blue, those with a SA circulation (if any) in red. 
	\label{fig:realTS}} 
\end{figure*}

The differences between the forecast temperature (forecast salinity) in the equatorial box and the polar box for the 100 ensemble members for the experiment without climate change are shown in \figRef{\ref{fig:realTS}}. It can be seen that during the DA period (first 18 years) both temperature and salinity differences follow the trend in the observations. However, they fail to reproduce the multiyear oscillation that is visible in the salinity and especially in the temperature observations. This is because such multiyear oscillations are missing in the surface forcing. This results in the model overestimating the salinity difference between equatorial and polar box thus underestimating the transport between the boxes. As explained in \secRef{\ref{sec:stommel_model}}, the autonomous solutions converge to one of at most two possible stable equilibria each associated with a different circulation regime. In the non-autonomous system, these equilibrium values will oscillate in time and the solution always ``chases'' one of them or switches from chasing one to chasing the other. The figure shows that in the absence of climate change the circulation remains in the TH regime for all ensemble members until at least 2104. Furthermore, the solutions tend to values that are close to current conditions: in the Stommel box model the solutions for the temperature differences lie above the difference at the end of the DA period. However, the differences are still within the range of contemporary temperature differences. The solutions for salinity differences are largersmaller at end of the DA period than in 2104, partially offsetting the decreasing trend over the DA period. The relative spread at the end of the century is also larger for the salinity differences than for the temperature differences. This can be explained by our observation in \secRef{\ref{sec:twin}} that the spread in salinity diffusion parameter is larger than that in the temperature diffusivity parameter. The values depend on the model parameters and consequently the uncertainty for future ocean salinities is relatively larger than that in ocean temperatures. 

\begin{figure*}[tbp]
	\centering
	\includegraphics[width=0.8\linewidth]{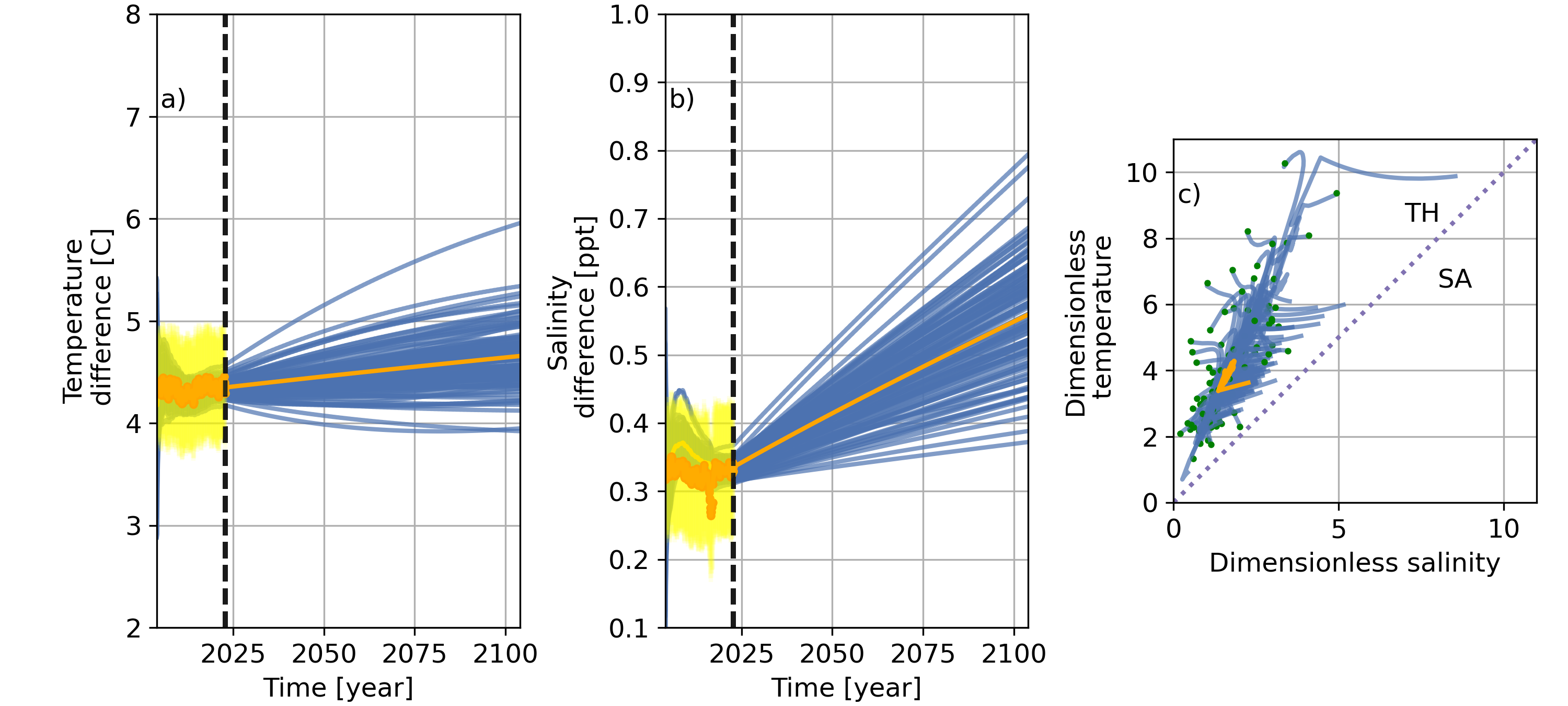}
	\caption{\raggedright As \figRef{\ref{fig:realTS}}, read left to right. Forecast time series starting in 2004 representing the equator-pole difference in a) temperature and b) salinity for the different ensemble members as well as c) the ensemble trajectories in a dimensionless temperature vs salinity phase portrait, but now climate change forcing is applied outside the assimilation period. The orange dots with bright yellow errorbars represent the observed values and their standard deviation. Solutions with a TH circulation are shown in blue, those with a SA circulation (if any) in red. End of the DA period is marked by a vertical dashed line \label{fig:realTSclima}}
\end{figure*}

The solutions in which after the DA period the polar surface temperature increases with $0.06 \,\unit{\celsius y^{-1}}$, the equatorial surface temperature with $0.03 \,\unit{\celsius y^{-1}}$ and in which the Green land ice sheet melts with a constant rate over 10 thousand years are shown in \figRef{\ref{fig:realTSclima}}. Compared to \figRef{\ref{fig:realTS}b} the salinity differences increase after the DA period (\figRef{\ref{fig:realTSclima}b}). This is due to the ice melt. Next to this, the temperature difference increase is smaller due to the differential increase in the surface temperature between the two boxes (\figRef{\ref{fig:realTSclima}}a). Apart from these points, there are no qualitative differences between the figures. In particular, none of the ensemble members changes circulation regime. 

\begin{figure}[tbp]
	\centering
	\includegraphics[width=0.6\textwidth]{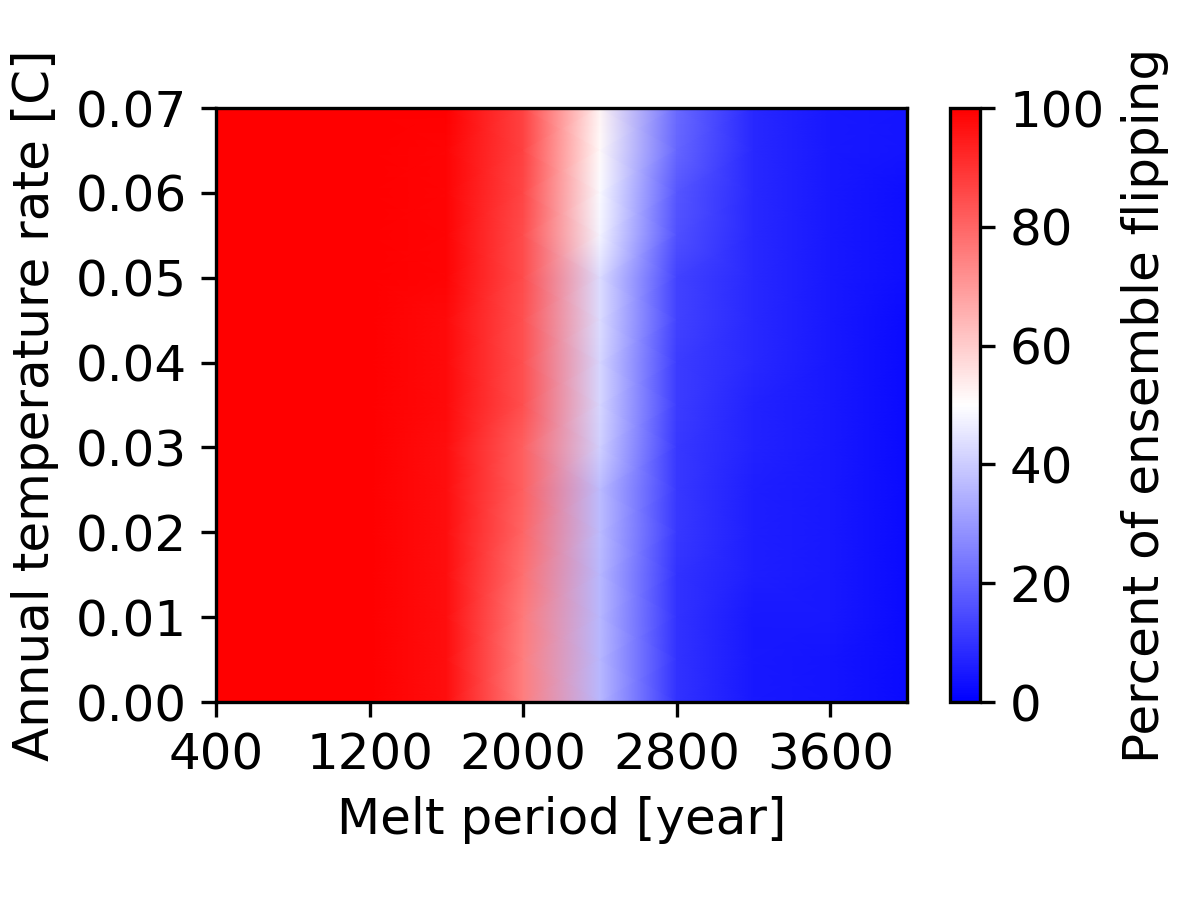}
	\caption{\raggedright Percentage of ensemble members that undergo circulation regime change. All runs use a similar setup to the experiment with climate change and real data assimilation, but using different ice melt periods (x-axis) and different warmings rates for the equator box (y-axis). The warming rate in the polar box is twice that of equator box.   Bright red denotes mostly flipping, while blue represents mostly remaining in the starting TH regime.}
	\label{fig:melt_clima_probs}
\end{figure}

The actual future rate of surface warming and ice melt will depend on, among other things, future carbon emission and are consequently uncertain. Therefore the DA experiment has been repeated several time for a range of plausible ice melting and between-box surface warming rate differences \cite{khan_greenland_2015,mouginot_forty-six_2019,mankoff_greenland_2020,you_warming_2021,xie_polar_2022}. The percentage of ensemble members in which the circulation regime changes from TH to SA are shown in \figRef{\ref{fig:melt_clima_probs}}. The figure shows that after DA the uncertainty in the regime is small: either all ensemble members are in the TH regime or in the SA regime. Next to this, the close-to-vertical slope of the boundary between the two regimes indicates that faster ice melting is more influential in pushing the solutions to the SA regime than faster warming over the pole than equator. Measured melt rates over the last decade vary between 286 to 487 $\unit{Gt y^{-1}}$ which corresponds to a melt period of 7-12 thousand years \cite{khan_greenland_2015,mouginot_forty-six_2019,mankoff_greenland_2020}. This would place the future climate firmly within the zone of \figRef{\ref{fig:melt_clima_probs}} in which circulation the AMOC does not reverse direction. 

\begin{figure}[h]
	\centering
	\includegraphics[width=3.5in]{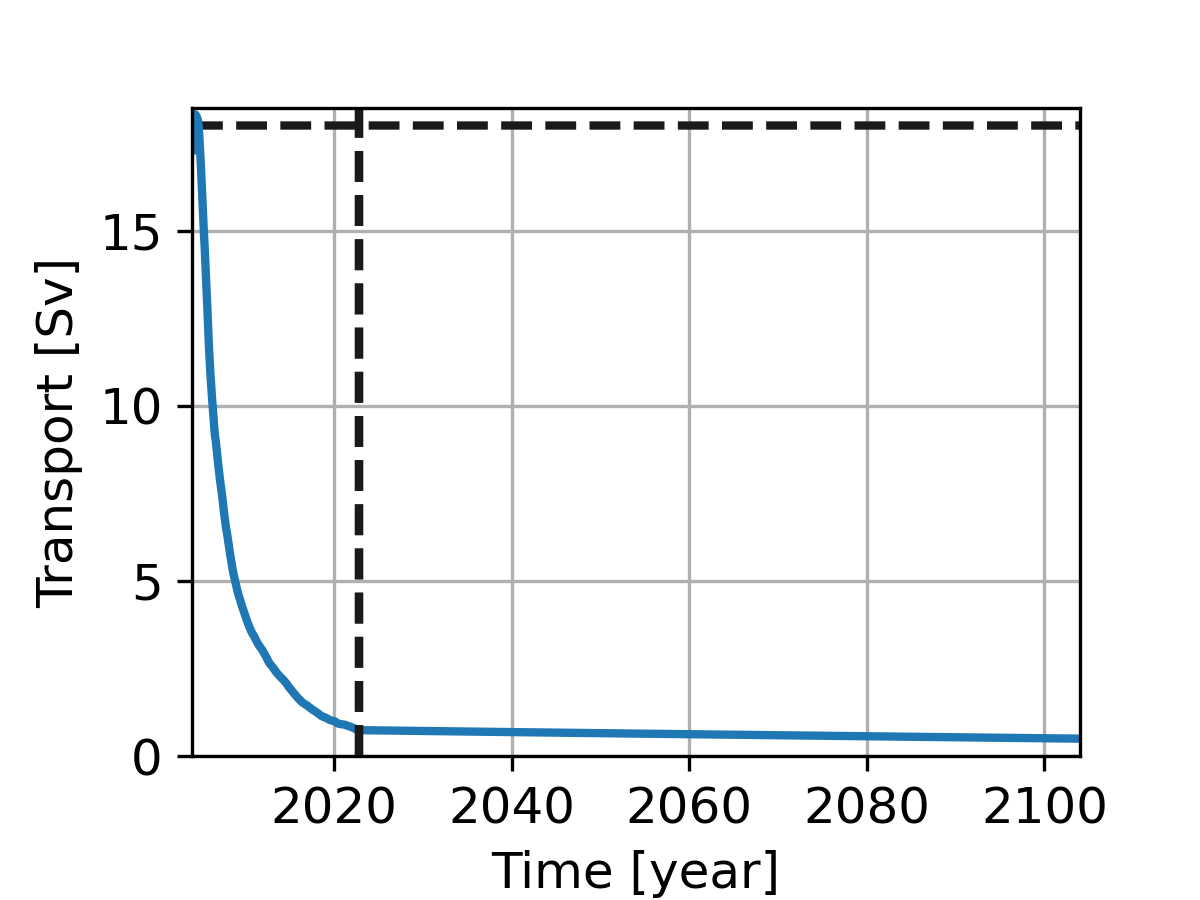}
	\caption{\raggedright The value of the between-box transport as function of time for the experiment with climate change. The vertical dashed line marks the end of the assimilation period. The horizontal dashed line indicates observed current transport \cite{MCC}. \label{fig:transport}}
\end{figure}

The transport between the two boxes is shown in \figRef{\ref{fig:transport}}. Transports are an order of magnitude smaller than the typical observational transect estimates over this period that lie in the range $10--25 \,\unit{Sv}$ \cite{danabasoglu_revisiting_2021}. The decrease is however not the result of climate change, but occurs during the DA period. I.e. the DA is not capable of simultaneously correcting temperatures, salinities and transports. The source of this discrepancy might lie in our box definition. This will be discussed in more detail in \secRef{\ref{sec:concl}}. The underestimation of current transports implies that the model is closer to the tipping point than the real AMOC. Consequently, the ice melt and warming rates in \figRef{\ref{fig:melt_clima_probs}} at which the circulation flips are in all likelihood an overestimation. 

\section{Conclusion and Discussions}\label{sec:concl}

We have applied the Stommel 2-box model to replicate the Atlantic meridional overturning circulation and investigated its solutions, the feasibility of retrieving the temperature diffusivity, salinity diffusivity and advection parameter by assimilating ocean temperature and salinity parameters and to forecast the likelihood of circulation regime change in the 21st century. 

We introduced seasonal variability and two of the effects of climate change into the model: warming of the ocean surface at different rate for the pole and the equator and the melt of the Greenland ice sheet.  Qualitatively, the Stommel box model \eqref{eq:dimensionless} mirrors certain aspects of the non-autonomous Stommel box model \eqref{eq:nonaut} close to the NSF (by setting $A=B=0$ in the non-autonomous case). The former presents a relatively generic and simple model with a region of bistability of two stable states that lose stability via a smooth saddle-node bifurcation while the latter undergoes an NSF bifurcation. {The non-autonomous Stommel box model presented provides a framework in which we can develop critical (tipping) points by introducing high frequency oscillatory forcing.} 

Using the setup, we tested whether it would be possible to reconstruct the true value of the model parameters by assimilating ocean temperature and salinity observation and use those parameters to forecast the future circulation regime of the AMOC. Our twin experiments showed that the DA method is indeed capable of reconstructing temperature diffusivity and the advection parameters from observations. It also bring the salinity diffusivity closer to the truth, but the 18 year DA period used is too short to convergence to the minimal obtainable error. Using observed temperature and salinity from the EN4 dataset, the most-likely values of the model parameters have been determined with the DA system indicating uncertainties of the same order of magnitude. In the absence of climate change future temperature and salinity is expected to lie in the same range of temperatures and salinities that have already been observed in the period 2004-2022. With the effects of climate, salinity differences between equatorial and pole box can move outside the current range because of the freshening caused by the melting Greenland ice sheet. Running the model with different climate change scenarios showed that this freshening is relatively more important to the direction of the AMOC circulation than changes in surface temperature. In particular, with the currently observed ice sheet melting rates, the model does not expect a reversal of the circulation direction to occur. 

Though an interesting proof-of-concept that data assimilation can help with reconstructing Stommel model parameters, the forecasts produced in this study should be interpreted with some care. The simplicity of the Stommel box model means that it cannot represent the wide range of processes occurring in the North-Atlantic, nor can it represent the North-Atlantic geometry. In particular, we have seen that our model underestimates the transports. {In other box-based studies \cite{rahmstorf_freshwater_1996,drijfhout_stability_2011,wood_observable_2019, chapman_quantifying_2024} the question of the transport is avoided by fitting temperature and salinity in boxes to a GCM with realistic levels for the AMOC transport. Consequently, defiencies in the dynamics of the Stommel model do not impact the parameter fit. In other multibox studies comparison with observed values is not even attempted \cite{monahan_stabilization_2002}. In \cite{gnanadesikan_simple_1999,kuhlbrodt_driving_2007} it has been argued that not the buoyancy forcing, but wind stresses are the primary driver behind the AMOC. And indeed \cite{cimatoribus_meridional_2014} has shown that the transport in a box-model increases strongly with increasing wind stresses. The fact such wind forcing is entirely absent in the Stommel model might explain why our transport estimates are an order of magnitude smaller than those observed. The omission of the southern Atlantic ocean in our setup might be another contributing factor to the errors in AMOC transports. Our two box model considers the North-Atlantic to be a closed system and does not account for the southward transport of North-Atlantic Deep Water past the equator towards the Antarctic Circumpolar Current nor for the northward flow of Antarctic Bottom Water and waters from the Angulhas current \cite[sec. 14.2.14]{talley_descriptive_2011} which can increase transports by creating diapycnal mixing \cite{cimoli_significance_2023}. Accommodating these concepts requires more complex conceptual models like the Rooth model \cite{rooth_hydrology_1982} or full climate models. Finally, we noticed that the volume of the boxes  can have profound impact on the transports and this dependence deserves more attention in future box-model studies.}

As a last note, the simplicity of the chosen reduced order model makes our analysis possible as we are able to initialize large ensembles at low computational costs. This setup can serve as a baseline in complexity with more detailed versions of the model available for further studies on higher order dynamics.

\appendix
\section{Lack of Periodic Solutions in the Stommel box model}
In this appendix, we go over why we can rule out the existence of periodic solutions in the autonomous Stommel box model. Doing so guarantees that in the non-autonomous model, when we do the twin experiments and begin assimilating real data, the only oscillatory behavior will come from the periodic forcing. The Poincaré-Bendixson theorem asserts that in a continuous dynamical system, any confined solution will inevitably approach either a fixed point or a limit cycle. In the next few paragraphs we will show that we get the former, not the latter in \eq{\eqref{eq:dimensionless}}, where we take $0<\eta_3<1$, $\eta_1>\frac{\eta_3}{1-\eta_3}$, with $T,S\geq0$.

To start, we have that
\begin{equation*}
	\begin{split}
		\frac{dT}{dt} =& f(T,S) = \eta_1 - T - T|T-S|\\
		\frac{dS}{dt} =& g(T,S) = \eta_2 - \eta_3 S - S|T-S| 
	\end{split}
\end{equation*}

We define our Dulac function $\phi(T,S)=1$. Applying Dulac's criterion, we get:
\begin{equation}
	\frac{\partial (\phi f)}{\partial T} + \frac{\partial (\phi g)}{\partial S} = \frac{\partial f}{\partial T} + \frac{\partial g}{\partial S}
	= -1 - \eta_3 - 3|T-S|
	\label{eq:dfdg}
\end{equation}
This is always strictly negative, and so by the Bendixson-Dulac theorem, we get that there cannot be any closed trajectories either entirely where $T>S>0$, as well as in the area where $S>T>0$. However, there could be a closed trajectory which crosses between the two areas, over the line $T=S$.

Now, assume we did have such a trajectory, a periodic solution $C$ which enclosed a region $D$ within $T,S\geq0$, and let us call the part part of $D$ where $T>S$, $D_1$, and the part where $S<T$ as $D_2$, so that $\partial D_1 = C_1 \cup C_3$ and $\partial D_2 = C_2 \cup (-C_3)$. 
\\
As \eq{\ref{eq:dfdg}} allows for a continuous extension from $D_1$ ($D_2$) to its compact closure $\overline{D}_1$ ($\overline{D}_2$), $\int\int_{\overline{D}_1}(\frac{\partial f}{\partial T}+\frac{\partial g}{\partial S})\:dT\:dS$ ($\int\int_{\overline{D}_2}(\frac{\partial f}{\partial T}+\frac{\partial g}{\partial S})\:dT\:dS$) exists and is, as shown above, strictly smaller than zero. Applying Green's Theorem now gives:
\begin{equation*}
	\footnotesize
	\begin{aligned}
		 \iint_D \left(\frac{\partial f}{\partial T} + \frac{\partial g}{\partial S}\right) dT dS 
		= & \iint_{D_1} \left(\frac{\partial f}{\partial T} + \frac{\partial g}{\partial S}\right) dT dS 
		 + \iint_{D_2} \left(\frac{\partial f}{\partial T} + \frac{\partial g}{\partial S}\right) dT dS \\
		&= \int_{C_1\cup C_3} (-g dT + f dS) + \int_{C_2\cup(-C_3)} (-g dT + f dS) \\
		&= \int_{C_1} (-g dT + f dS) + \int_{C_3} (-g dT + f dS) 
		 + \int_{C_2} (-g dT + f dS) + \int_{(-C_3)} (-g dT + f dS)
	\end{aligned}
\end{equation*}

But rearranging this, the second and fourth terms simply cancel out. Furthermore, by definition, we have that $\frac{dT}{dt}=f(T,S)$ and $\frac{dS}{dt}=g(T,S)$, so we get:
\begin{equation*}
	\footnotesize
	\begin{aligned}
		\int\int_{D}(\frac{\partial f}{\partial T}+\frac{\partial g}{\partial S})\:dT\:dS=&\int_{C_1}(-g\:dT+f\:dS) + \int_{C_2}(-g\:dT+f\:dS)\\
		=&\int_C (-g\:dT+f\:dS)\\
		=&\int_C (-g\frac{dT}{dt}+f\frac{dS}{dt})\:dt\\
		=&\int_C (-\frac{dS}{dt}\frac{dT}{dt}+\frac{dT}{dt}\frac{dS}{dt})\:dt\\
		=&\int_C 0\:dt\\
		=& 0
	\end{aligned}
\end{equation*}
But earlier we showed this was less than zero, so we have a contradiction, so no such loop $C$ can exist.

Thus, in the autonomous Stommel box model, it is impossible for a periodic solution to exist for $T,S \geq 0$.

\section{Ensemble Kalman filter}
\label{app:kalman}

Let 
\begin{align*}
          \Vec{x}_{t}=& [T_{e*}(t),T_{p*}(t),S_{e*}(t),S_{p*}(t),\log \frac{\kappa_{T}}{\kappa_{ref}}, \log \frac{\kappa_{S}}{\kappa_{ref}}, \log \frac{\gamma}{\kappa_{ref}}]^{\T}  
\end{align*}
\noindent {represent a possible state of the system at time $t$ with $\kappa_{ref}=1 \, \unit{m^{2}s^{-1}}$ a reference values to nondimensionlize the diffusion and advection coefficients. As this state does not only include the state variables ($T_{e*}(t)$,$T_{p*}(t)$, $S_{e*}(t)$, $S_{p*}(t)$) but also the model parameters $\kappa_{T}$, $\kappa_{S}$ and $\gamma$, it is referred to as the augmented state. Let $\int_{\mathcal{X}} p(\Vec{x}_{0}) \,\mathrm{d}\Vec{x}_{0}$ be the probability that the true state of the system at the beginning of the model run lies in the set $\mathcal{X}$. Using Bayes' theorem and the chain rule for probabilities it can be shown \cite{evensen_data_2022-1} that 
\begin{equation}
	\underbrace{p(\Vec{x}_{t}|\Vec{d}_{0:t})}_{posteriori} \sim p(\Vec{d}_{t}|\Vec{x}_{t})\underbrace{\int p(\Vec{x}_{t}|\Vec{x}_{t-1})p(\Vec{x}_{t-1}|\Vec{d}_{0:t-1}) \,\mathrm{d}\Vec{x}_{t-1}}_{apriori}
	\label{eq:bayes}
\end{equation}
where $\sim$ indicates equivalence up to a normalisation constant, $d_{t}$ a vector with observations at timestep $t$ and $d_{0:t}$ a vector of all observations upto and including timestep $t$. Once the apriori is specified for the initial timestep $t=0$, equation~\ref{eq:bayes} posteriori and apriori distributions for $t>0$ can be constructed by induction using Equation~\ref{eq:bayes}. This occurs in two steps. The {\it prediction step} in which the posteriori distribution at time $t-1$ is transformed into the apriori distribution at time $t$, and the {\it update step} in which the DA correction is applied to transform the apriori distribution at time $t$ into the posteriori distribution at the same time. 

In an ensemble Kalman filter it is assumed that the apriori and posteriori distributions can be approximated by a Gaussian with its mean $\Vec{\mu}_{x,t}=\frac{1}{M} \sum_{m=1}^{M} \Vec{x}_{t}^{(m)}$ and covariance $\Sigma_{x}=\frac{1}{M-1}\sum_{m=1}^{M}(\Vec{x}^{(m)}_{t}-\Vec{\mu}_{x,t})(\Vec{x}^{(m)}_{t}-\Vec{\mu}_{t})^{\T}$ estimated from an $M$-member ensemble of model runs $\{\Vec{x}_{t}^{(m)}:1\leq m \leq M\}$. Each ensemble member represents an equally likely realisation of these probability distributions. The integral in the apriori distribution in Equation~\ref{eq:bayes} is approximated using a Monte-Carlo method. For each ensemble $\Vec{x}^{(n)}_{t-1}$ an ensemble member $\Vec{x}^{(m)}_{t}$ is sampled according to $p(\Vec{x}_{t}|\Vec{x}_{t-1})$. In this work the model is assumed to be without model error and hence $p(\Vec{x}_{t}|\Vec{x}_{t-1})=\delta(\Vec{x}_{t}-\mathcal{M}_{t-1 \to t}\Vec{x}_{t-1})$ with $\delta$ the delta function and $\mathcal{M}_{t\to t-1}$ the model operator that propagates the ocean temperatures and salinities from time $t-1$ to $t$ according to equations \ref{eq:main}-\ref{eq:eos}. $\mathcal{M}_{t-1 \to t}$ whilst keeping the model parameters constant. I.e. model parameters are only changed during the {\it update step}.

The probability density distribution for the observations $p(\Vec{d}_{t}|\Vec{x}_{t})$ is assumed to be a Gaussian with mean $\mu_{d}=H \Vec{\mu}_{t,x}$ and covariance $\Sigma_{d}$. In our case the observation operator $H$ is such that $H \Vec{x}_{t}=[T_{p*}(t),T_{e*}(t),S_{p*}(t),S_{e*}(t)]^{\T}$ and $\Sigma_{d}$ is assumed to be diagonal with its diagonal elements (variances) given by \eq{\ref{eq:T_sigma}}. 

It can be shown that under these assumptions, the posteriori distribution is again a Gaussian with mean $\Vec{\mu}_{\hat{x},t}=\frac{1}{M} \sum_{m=1}^{M} \Vec{\hat{x}}^{(m)}_{t}$ and covariance $\Sigma_{\hat{x},t}=\frac{1}{M-1} (\Vec{\hat{x}}_{t}^{(m)}-\Vec{\mu}_{\hat{x},t})(\Vec{\hat{x}}_{t}^{(m)}-\Vec{\mu}_{\hat{x},t})^{\T}$
and
\begin{equation}
	\Vec{\hat{x}}_{t}^{(n)} = \Vec{x}_{t}^{(n)}+ K(\Vec{d}_{t} - H\Vec{\mu}_{t,x}) + A Q \Lambda^{\frac{1}{2}} Q^{\T}\Vec{\hat{e}}^{(n)},
\end{equation}
with $\Vec{\hat{e}}^{(n)}$ the $n$th unit vector, $A$ a matrix having $\Vec{x}_{t}^{(n)}-\Vec{\mu}_{t,x}$ as its $n$th column and $Q \Lambda Q^{\T}$ the singular value decomposition of  $I - A^{\T} H^{\T} ((m-1)\Sigma_{d}+H A A^{\T} H^{\T})^{-1} H A$ \cite{evensen_data_2022-1} and $K=AA^{\T}(H A A^{\T} H^{\T}+(M-1)\Sigma_{d})^{-1}$ the Kalman gain. 

The most-likely state, or mode, of the apriori and posteriori distributions are referred to as the forecast and analysis respectively. As the ocean temperatures and salinities are assumed to follow Gaussian distributions, their forecasts and analysis coincide with their ensemble mean. For the advection and diffusion parameters, however, we corrected their logarithms instead of the values themselves. This so-called log-transform ensures that the parameters themselves remain strictly positive. For a log-normal distribution, the mode is not given by $e^{\Vec{\mu}_{p}}$ as naively might be expected, but by $e^{\Vec{\mu}_{p}-\Sigma_{p} \Vec{1}}$ with $\Vec{\mu}_{p}$ the 3-dimensional vector containing the ensemble mean of the logarithms of the model parameters and $\Sigma_{p}$ the $3\times3$ matrix with their ensemble covariance~\cite{fletcher_mixed_2010}. I.e. $\mu_{p}$ and $\Sigma_{p}$ are the last 3-elements or $3\times3$-block of the $\mu_{x}$ and $\Sigma_{x}$ respectively.}

\section*{Data Availability}
The MetOffice EN4 objective analysis temperature and salinity data with XBT and MBT corrections used to create observations for assimilation as well as the surface forcing is available at \url{https://hadleyserver.metoffice.gov.uk/en4/download-en4-2-2.html} as the \emph{EN.4.2.2.analyses.g10.XXXX.zip} files. The \emph{hadley\_obs.py} script used to process these observations can be found in the \emph{inversion} branch of the \cite{pasmans_dapperinversion_2023} GIT repository. Also included in this fork of DAPPER \cite{raanes_dapper_2023} are the code for the Stommel model (\emph{dapper/mods/stommel/\_\_init\_\_.py}) as well as the scripts used to carry out the data assimilation (\emph{dapper/mods/stommel/experiment.py}) using these observations and the scripts for the twin experiments. 

\section*{Acknowledgement}
We would like to thank John Gemmer and Christopher K.R.T. Jones for starting and funding this research project during the 2023 MCRN Summer School and Research Program under NSF grant DMS-1722578. Ivo Pasmans would like to acknowledge the support of the Scale-Aware Sea Ice Project (SASIP) funded by Schmidt Sciences – a philanthropic initiative that seeks to improve societal outcomes through the development of emerging science and technologies - grant G-24-66154. Emmanuel Fleurantin was supported by NSF grant DMS-2137947. A part of this project originated from the Mason Experimental Geometry Lab (MEGL) at George Mason University. A special thanks goes out to those MEGL students who started with the twin experiment setup in the Fall 2022 semester.

\bibliographystyle{unsrt}  
\bibliography{Bibliography,Stommel_pasmans}


%
%

%


\end{document}